\RequirePackage{ifpdf}
\ifpdf 
\documentclass[pdftex]{sigma}
\else
\documentclass{sigma}
\fi

\newcommand{\Z}{{\Bbb Z}} 
\newcommand{\C}{{\Bbb C}} 

\newcommand{\cA}{{\cal A}}
\newcommand{\cU}{{\cal U}}
\newcommand{\cZ}{{\cal Z}}

\newcommand{\cR}{{\cal R}}

\newcommand{\ve}{\varepsilon}

\newcommand{\la}{\lambda}

\newcommand{\De}{{\Delta}}

\newcommand{\Hom}{\operatorname{Hom}}

\newcommand{\nn}{{\nonumber}}
\newcommand{\bea}{\begin{eqnarray}}
\newcommand{\ena}{\end{eqnarray}}
\newcommand{\beit}{\begin{itemize}}
\newcommand{\enit}{\end{itemize}}

\newcommand{\be}{\begin{eqnarray*}}
\newcommand{\en}{\end{eqnarray*}}
\newcommand{\lb}[1]{\label{#1}}

\newcommand{\ds}[1]{{\displaystyle #1 }}

\newcommand{\kuru}{\curvearrowleft} 
\newcommand{\id}{{\rm id}}
\newcommand{\Ad}{{\rm Ad}}

\newcommand{\wt}{{\rm wt}}

\newcommand{\tp}{{\tilde{p}}}

\newcommand{\BW}[5]
{\left(\begin{array}{cc}#1 & #2 \cr #3 & #4 \cr\end{array} \Biggl|
#5\right)}

\def\infq4p#1{{(#1;q^4,p)_\infty}}

\newcommand{\tPsi}{\widetilde{\Psi}}

\newcommand{\al}{\alpha}
\newcommand{\e}{\epsilon}
\newcommand{\vep}{\varepsilon}

\newcommand{\bLa}{\bar{\Lambda}}

\newcommand{\mmatrix}[1]{\begin{matrix} #1 \end{matrix}}
\newcommand{\mat}[1]{\left(\mmatrix{#1}\right)}

\font\teneufm=eufm10 \font\seveneufm=eufm7 \font\fiveeufm=eufm5
\newfam\eufmfam
\textfont\eufmfam=\teneufm \scriptfont\eufmfam=\seveneufm
\scriptscriptfont\eufmfam=\fiveeufm

\let\goth\frak
\newcommand{\slth}{\widehat{\goth{sl}}_2}
\newcommand{\slt}{\goth{sl}_2}
\newcommand{\slnh}{\widehat{\goth{sl}}_N}

\newcommand{\g}{\goth{g}}
\newcommand{\gb}{\goth{b}}

\newcommand{\Aqp}{{\cal A}_{q,p}}

\newcommand{\Bqla}{{{\cal B}_{q,\lambda}}}

\newcommand{\h}{\goth{h}}

\font\fourteeneufm=eufm10 scaled\magstep2    
\font\seventeeneufm=eufm10 scaled\magstep3   

\newcommand{\gbig}{\mbox{\fourteeneufm g}} 
\newcommand{\gBig}{\mbox{\seventeeneufm g}} 
\makeatletter \@addtoreset{equation}{section} \makeatother

\newtheorem{thm}{Theorem}[section]
\newtheorem{prop}[thm]{Proposition}
\newtheorem{lem}[thm]{Lemma}

\newtheorem{conj}[thm]{Conjecture}

\newtheorem{df}{Definition}[section]
\newtheorem{dfn}[thm]{Definition}

\begin{document}

\allowdisplaybreaks

\renewcommand{\PaperNumber}{091}

\FirstPageHeading

\renewcommand{\thefootnote}{$\star$}

\ShortArticleName{Dynamical $R$ Matrices of Elliptic Quantum
Groups}

\ArticleName{Dynamical $\boldsymbol{R}$ Matrices of Elliptic Quantum Groups  \\
  and Connection Matrices for the $\boldsymbol{q}$-KZ Equations\footnote{This paper is a contribution to the Proceedings of
the O'Raifeartaigh Symposium on Non-Perturbative and Symmetry
Methods in Field Theory
 (June 22--24, 2006, Budapest, Hungary).
The full collection is available at
\href{http://www.emis.de/journals/SIGMA/LOR2006.html}{http://www.emis.de/journals/SIGMA/LOR2006.html}}}

\Author{Hitoshi KONNO} 
\AuthorNameForHeading{H.~Konno}

\Address{Department of Mathematics, Graduate School of Science, \\
Hiroshima University, Higashi-Hiroshima 739-8521, Japan}
\Email{\href{mailto:konno@mis.hiroshima-u.ac.jp}{konno@mis.hiroshima-u.ac.jp}}

\ArticleDates{Received October 02, 2006, in f\/inal form November
28, 2006; Published online December 19, 2006}

\Abstract{For any af\/f\/ine Lie algebra $\g$, we show that any
f\/inite dimensional representation of the universal dynamical $R$
matrix  $\cR(\la)$ of the elliptic quantum group ${\cal
B}_{q,\lambda}(\g)$ coincides with a corresponding connection
matrix for the solutions of the $q$-KZ equation associated with
$U_q(\g)$. This provides a general connection between ${\cal
B}_{q,\lambda}(\g)$ and the elliptic face (IRF or SOS) models. In
particular, we construct  vector representations of $\cR(\la)$ for
$\g=A_n^{(1)}$, $B_n^{(1)}$, $C_n^{(1)}$, $D_n^{(1)}$, and show
that they coincide with the face weights derived by Jimbo, Miwa
and Okado. We hence conf\/irm the conjecture by Frenkel and
Reshetikhin.}


\Classification{33D15; 81R50; 82B23}

\section{Introduction}
The quantum group $U_q(\g)$ is one of the
 fundamental structures appearing in the wide class of trigonometric quantum
integrable systems. Among others, we remark the following two
facts.
\begin{itemize}\itemsep=0pt
\item[1)] For $\g$ being af\/f\/ine Lie algebra, f\/inite
dimensional representations of $U_q(\g)$ allow a syste\-matic
derivation of trigonometric solutions of the Yang--Baxter equation
(YBE) \cite{Jimbo,Jimbo86}. \item[2)] A combined use of f\/inite
and inf\/inite dimensional representations allows us to formulate
trigonometric vertex models
 and  calculate correlation functions~\cite{JM}.
\end{itemize}

To extend this success to elliptic systems is our basic aim. In
this paper, we consider a~problem analogous to 1). As for
developments in the direction 2), we refer the reader to the
papers~\cite{Konno, JKOS2, KK03, KK04, KKW}.  We are especially
interested in the two dimensional exactly solvable lattice models.
There are two types of elliptic solvable lattice models. The
vertex type and the face (IRF or SOS) type. The vertex type
elliptic solutions to the YBE were found by Baxter~\cite{Baxter}
and Belavin~\cite{Belavin}. These are classif\/ied as the elliptic
$R$ matrices of the type $A_n^{(1)}$. The face type elliptic
Boltzmann weights associated with $A_1^{(1)}$ were f\/irst
constructed by Andrews--Baxter--Forrester~\cite{ABF}, and extended
to $A_n^{(1)}$, $B_n^{(1)}$, $C^{(1)}_n$, $D_n^{(1)}$ by
 Jimbo--Miwa--Okado~\cite{JMOAtype, JMO},
 to $A^{(2)}_{2n}$, $A^{(2)}_{2n-1}$ by Kuniba~\cite{Kuniba},
and to~$G^{(1)}_2$ by Kuniba--Suzuki~\cite{KS}.

Concerning the elliptic face weights, Frenkel and Reshetikhin made
an interesting observation~\cite{FR}  that the connection matrices
for the solution of the $q$-KZ equation associated with $U_q(\g)$
($\g$: af\/f\/ine Lie algebra) provide elliptic solutions to the
face type YBE. They also conjectured that the connection matrices
in
 the vector representation are equal to Jimbo--Miwa--Okado's face weights for
$\g=A_n^{(1)}, B_n^{(1)}, C^{(1)}_n, D_n^{(1)}$. In order to
conf\/irm this conjecture, one needs to solve the $q$-KZ equation
of general level and f\/ind connection matrices. Within our
knowledge, no one has yet conf\/irmed it. Instead of doing this,
Date, Jimbo and Okado \cite{DJO} considered the face models
def\/ined by taking the connection matrices as  Boltzmann weights.
They showed that the one-point function of such models is given by
the branching function associated with $\g$. The same property of
the one-point function had been discovered in Jimbo--Miwa--Okado's
$A_n^{(1)}$, $B_n^{(1)}$, $D_n^{(1)}$ face models.

An attempt to formulate elliptic algebras  was f\/irst made by
Sklyanin~\cite{Sklyanin}. He considered an algebra def\/ined by
the $RLL$-relation associated with Baxter's elliptic $R$-matrix.
It was extended to the
 elliptic algebra $\Aqp(\slth)$ by Foda et al.~\cite{FIJKMY},
 based on a central extension of Sklyanin's $RLL$-relation.
In the same year, Felder proposed a face type elliptic algebra
$E_{\tau,\eta}(\g)$ associated with the dynamical
$RLL$-relation~\cite{Felder}. Jimbo--Miwa--Okado's  elliptic
solutions to the face type YBE were interpreted there as the
dynamical $R$ matrices. We classify the former elliptic algebra
the vertex-type and the latter the face-type. Another formulation
of the face type
 elliptic algebra was discovered by the  author~\cite{Konno}. It is based on an elliptic deformation of
the Drinfel'd currents.

A coalgebra structure of these elliptic algebras was clarif\/ied
in the works by Fr\o nsdal~\cite{Fronsdal},
Enriquez--Felder~\cite{EF} and
Jimbo--Konno--Odake--Shiraishi~\cite{JKOS}. It is based on an idea
of quasi-Hopf deformation~\cite{Dri90} by using the twistor
operators satisfying the shifted cocycle condition~\cite{BBB}. In
this formulation, we regard the coalgebra structures of the vertex
and the face type elliptic algebras as two dif\/ferent quasi-Hopf
deformation of the corresponding af\/f\/ine quantum
group~$U_q(\g)$. We call the resultant quasi-Hopf algebras the
elliptic quantum groups of the vertex type $\Aqp(\slnh)$ and the
face type $\Bqla(\g)$. A detailed description for the face type
case is reviewed in Section~2. For the face type, a dif\/ferent
coalgebra structure as a $\h$-Hopf algebroid was developed by
Felder, Etingof and Varchenko \cite{FV, EV1, EV2},
Koelink--van~Norden--Rosengren~\cite{KNR}.

One of the advantages of the quasi-Hopf formulation is that it
allows
 a natural derivation of the universal dynamical $R$ matrix from one of
$U_q(\g)$ as a twist. However, a disadvantage is that there are no
a priori reasons for the resultant universal $R$ matrix to yield
elliptic $R$ matrices. One needs to check this point in all
representations. We have done this for the vector representations
of $\Aqp(\slth)$ and  $\Bqla(\slth)$, which led to
 Baxter's elliptic $R$ matrix and
Andrews--Baxter--Forrester's elliptic face weights,
respectively~\cite{JKOS}. The same checks for the face weights
were also done in the cases $\g=A_n^{(1)},
A_2^{(2)}$~\cite{KK03,KK04}.

The aim of this paper is to overcome this disadvantage by
clarifying
 the following point concerning the face type.
\begin{itemize}\itemsep=0pt
\item[i)] Any representations of the universal dynamical $R$
matrix of $\Bqla(\g)$ are equivalent to the corresponding
connection matrices for the $q$-KZ equation of $U_q(\g)$.
\end{itemize}
The connection matrices are known to be elliptic. See Theorem
\ref{connectMat}. In addition, we show
\begin{itemize}\itemsep=0pt
\item[ii)] For $\g=A_n^{(1)}, B_n^{(1)}, C^{(1)}_n, D_n^{(1)}$,
the vector representation of the the universal dynamical $R$
matrix of $\Bqla(\g)$ is equivalent to Jimbo--Miwa--Okado's
elliptic face weight up to a gauge transformation.
\end{itemize}

Combining i) and ii), we conf\/irm Frenkel--Reshetikhin's
conjecture on the equivalence between the connection matrices and
Jimbo--Miwa--Okado's face weights. For the purpose of showing i),
we follow the idea by Etingof and Varchenko~\cite{EV2}, and give
an exact relation between the face type twistors and the highest
 to highest  expectation values of the composed vertex
operators (fusion matrices) of $U_q(\g)$. To show ii), we solve
the dif\/ference equation for the face type twistor, which is
equivalent to the $q$-KZ equation of general level.

This paper is organized as follows. In the next section, we
summarize some basic facts on the af\/f\/ine quantum groups
$U_q(\g)$ and the face type elliptic quantum groups $\Bqla(\g)$.
In Section~3, we introduce the vertex operators of $U_q(\g)$ and
fusion matrices.  We discuss equivalence between the face type
twistors and the fusion matrices. Then in Section~4, we show
equivalence between the dynamical $R$ matrices of $\Bqla(\g)$ and
the connection matrices for the $q$-KZ equation of  $U_q(\g)$ in
general f\/inite dimensional representation. Section~5 is devoted
to a discussion on an equivalence between the vector
representations of the the universal dynamical $R$ matrix of
$\Bqla(\g)$ and Jimbo--Miwa--Okado's elliptic face weights for
$\g=A_n^{(1)}, B_n^{(1)}, C^{(1)}_n, D_n^{(1)}$.

\renewcommand{\thefootnote}{\arabic{footnote}}
\setcounter{footnote}{0}

\section[Affine quantum groups $U_q(\gBig)$ and elliptic quantum groups ${\cal B}_{q,\la}(\gBig)$]{Af\/f\/ine
quantum groups $\boldsymbol{U_q(\gBig)}$\\ and elliptic quantum
groups $\boldsymbol{{\cal B}_{q,\la}(\gBig)}$}

\subsection[Affine quantum groups $U_q(\gbig)$]{Af\/f\/ine quantum groups $\boldsymbol{U_q(\gbig)}$}

Let $\g$ be an af\/f\/ine Lie algebra associated with a
generalized Cartan
 matrix $A=(a_{ij})$, $i,j\in I=\{0,1,\ldots,n\}$.
  We f\/ix an invariant inner product $(\cdot |\cdot )$ on the
  Cartan subalgebra $\h$ and identify $\h^*$ with $\h$ through $(\cdot |\cdot )$.
We follow the notations and conventions in \cite{Kac} except for
$A_{2n}^{(2)}$, which we def\/ine in such a way that  the order of
the vertices of the Dynkin diagram is reversed from the one
in~\cite{Kac}. We hence  have $a_0=1$ for all $\g$.
Let $\{\alpha_i\}_{i\in I}$ be a set of simple roots and set
$h_i=\al_i^\vee$. We have  $a_{ij}=\langle\al_j,h_i\rangle
=\frac{2(\al_i| \al_j)}{(\al_i|\al_i)}$ and $d_ia_{ij}=a_{ij}d_j$
with $d_i=\frac{1}{2}(\al_i|\al_i)$. We denote the canonical
central element by $c=\sum_{i\in I} a_i^{\vee}h_i$ and the null
root by $\delta=\sum_{i\in I} a_i\alpha_i$. We set
\begin{gather*}
Q=\Z\al_0\oplus \cdots \oplus \Z\al_n,\qquad Q_+=\Z_{\geq 0}\al_0\oplus \cdots \oplus \Z_{\geq 0}\al_n,\\
P=\Z\Lambda_0\oplus \cdots \oplus \Z\Lambda_n \oplus \Z
\delta,\qquad P^*=\Z h_0\oplus \cdots \oplus \Z h_n \oplus \Z d
\end{gather*}
and impose the pairings
\[
\langle \alpha_i,d\rangle=a_0^\vee\delta_{i,0},\qquad
\langle\Lambda_j,h_i\rangle =\frac{1}{a_0^\vee}\delta_{ij},\qquad
\langle\Lambda_j,d\rangle=0,
 \qquad (i,j\in I).
\]
The $\Lambda_j$ are the fundamental weights. We also use
$P_{cl}=P/\Z\delta$, $(P_{cl})^*=\oplus_{i=0}^n\Z h_i \subset
P^*$. Let $cl:P\to P_{cl}$  denote the canonical map and def\/ine
$af:P_{cl}\to P$ by $af(cl(\al_i))=\al_i \ (i\not=0)$ and
$af(cl(\Lambda_0))=\Lambda_0$ so that $cl\circ af=\id$ and
$af(cl(\al_0))=\al_0-\delta$.

\begin{dfn}
The quantum affine algebra $U_q=U_q(\g)$ is an associative algebra
over $\C(q^{1/2})$ with $1$ generated by the elements $e_i$, $f_i$
$(i\in I)$ and $q^h\ (h\in P^*)$ satisfying the following
relations
\begin{gather*}
q^0=1,\ q^hq^{h'}=q^{h+h'}\qquad (h, h'\in P^*),\\
q^he_i q^{-h}=q^{\langle \al_i,h\rangle}e_i,\qquad
q^hf_i q^{-h}=q^{-\langle \al_i,h\rangle}f_i,\\
e_if_j-f_je_i=\delta_{ij}\frac{t_i-t_i^{-1}}{q_i-q_i^{-1}},\\
\sum_{m=0}^{1-a_{ij}}(-1)^m \left[\mmatrix{1-a_{ij}\cr m}\right]_{q_i} e_i^{1-a_{ij}-m}e_je_i^m=0\qquad (i\not=j),\\
\sum_{m=0}^{1-a_{ij}}(-1)^m \left[\mmatrix{1-a_{ij}\cr
m}\right]_{q_i} f_i^{1-a_{ij}-m}f_jf_i^m=0\qquad (i\not=j).
\end{gather*}
Here $q_i=q^{d_i}$, $t_i=q_i^{h_i}$, and
\begin{gather*}
[n]_x=\frac{x^n-x^{-n}}{x-x^{-1}},\qquad
[n]_x!=[n]_x[n-1]_x\cdots[1]_x,\qquad
 \left[\mmatrix{n\cr m \cr}\right]_{x}=\frac{[n]_x!}{[m]_x![n-m]_x!}.
\end{gather*}
\end{dfn}

The algebra $U_q$ has a Hopf algebra structure with
comultiplication $\De$, counit $\epsilon$ and antipode~$S$
def\/ined by
\begin{gather}
\De(q^h)=q^h\otimes q^h,\nn\\
\De(e_i)=e_i\otimes 1+t_i\otimes e_i,\nn\\
\De(f_i)=f_i\otimes t_i^{-1}+1\otimes f_i,\lb{coalgstr}\\
\e(q^h)=1,\qquad \e(e_i)=\e(f_i)=0,\nn\\
S(q^h)=q^{-h},\qquad S(e_i)=-t_i^{-1}e_i,\qquad S(f_i)=-f_it_i.\nn
\end{gather}

$U_q$ is a quasi-triangular Hopf algebra with
 the universal $R$ matrix $\cR$ satisfying
\begin{gather*}
\De^{op}(x)=\cR\De(x)\cR^{-1}\qquad \forall\, x\in U_q,\\
(\De \otimes \id)\cR=\cR^{13}\cR^{23},\qquad 
(\id \otimes \De)\cR=\cR^{13}\cR^{12}.
\end{gather*}
Here $\De^{op}$ denotes the opposite comultiplication,
$\De^{op}=\sigma\circ \De$ with $\sigma$ being the f\/lip of the
tensor components; $\sigma(a\otimes b)=b\otimes a$.
\begin{prop}
\begin{gather}
\cR^{(12)}\cR^{(13)}\cR^{(23)} =
\cR^{(23)}\cR^{(13)}\cR^{(12)},\lb{YBE}\\
(\e\otimes \id)\cR=(\id \otimes \e)\cR=1,\nn\\
(S\otimes \id)\cR=(\id \otimes S^{-1})\cR=\cR^{-1},\nn\\
(S\otimes S)\cR=\cR.\nn
\end{gather}
\end{prop}

Let $\{h_l\}$ be a basis of $\h$ and $\{h^l\}$ be its dual basis.
We denote by $U^+$ (resp. $U^-$) the subalgebra of $U_q$ generated
by $e_i$ (resp. $f_i$) $i\in I$ and set
\begin{gather*}
U^+_\beta=\{ x \in U^+|\ q^hxq^{-h}=q^{\langle\beta,h\rangle}x\  (h\in \h) \} ,\\
U^-_{-\beta}=\{ x \in U^-|\ q^hxq^{-h}=q^{-\langle
\beta,h\rangle}x\  (h\in \h) \}
\end{gather*}
for $\beta \in Q_+$. The universal  $R$ matrix has the
form~\cite{Tanisaki}
\begin{gather}
\cR=q^{-T}\cR_0,\qquad T=\sum_{l}h_l\otimes h^l,\nn\\
\cR_0=\sum_{\beta\in Q_+}q^{(\beta,\beta)}\big(q^{-\beta}\otimes
q^\beta \big)
(\cR_0)_\beta=1-\sum_{i\in I}(q_i-q_i^{-1})e_it_i^{-1}\otimes t_if_i+\cdots,\nn\\
(\cR_0)_\beta=\sum_j u_{\beta,j}\otimes u^j_{-\beta} \in
U_\beta^+\otimes U^-_{-\beta}, \lb{strR}
\end{gather}
where $\{u_{\beta,j}\}$ and $\{u^j_{-\beta}\}$ are bases of
$U^+_{\beta}$ and $U^-_{-\beta}$, respectively. Note that $T$ is
the canonical element of $\h \otimes \h$ w.r.t $(\cdot |\cdot )$
and $(\cR_0)_\beta$ is the canonical element of $U_\beta^+\otimes
U^-_{-\beta}$ w.r.t a certain Hopf paring.

We write $U'_q=U'_q(\g)$ for the subalgebra of $U_q$ generated by
$e_i$, $f_i$ $(i\in I)$ and $h\in (P_{cl})^*$.
Let $(\pi_V,V)$ be a f\/inite dimensional module over $U'_q$. We
have the evaluation representation $(\pi_{V,z}, V_z)$ of $U_q$ by
$V_z=\C(q^{1/2})[z,z^{-1}]\otimes_{\C(q^{1/2})}V$ and
\begin{gather*}
\pi_{V,z}(e_i)(z^n\otimes v)=z^{\delta_{i0}+n}\otimes
\pi(e_i)v,\qquad
\pi_{V,z}(f_i)(z^n\otimes v)=z^{-\delta_{i0}+n}\otimes \pi(f_i)v,\\
\pi_{V,z}(t_i)(z^n\otimes v)=z^{n}\otimes \pi(t_i)v,\qquad
\pi_{V,z}(q^d)(z^n\otimes v)=(qz)^{n}\otimes v, \\
\wt(z^{n}\otimes v)=n\delta+af(\wt(v)),
\end{gather*}
where $n\in \Z$, and  $v\in V$ denotes a weight vector whose
weight is $\wt(v) $. We write $vz^{n}=z^{n}\otimes v$ $(n\in \Z)$.

For generic $\la\in \h^*$, let $M_\la$ denote the irreducible
Verma module with the highest weight~$\la$. We have the weight
space decomposition $\ds{M_\la=\bigoplus_{\nu\in
\la-Q_+}(M_\la)_\nu}$. We write $\wt(u)=\nu $ for $u\in
(M_\la)_\nu$.

\subsection[Elliptic quantum groups ${\cal B}_{q,\la}(\gbig)$]{Elliptic quantum
groups $\boldsymbol{{\cal B}_{q,\la}(\gbig)}$}

Let $\rho\in \h$ be an element satisfying $(\rho|\alpha_i)=d_i$
for all $i\in I$. For generic $\la\in \h$, let us consider an
automorphism of $U_q$ given by
\[
\varphi_\la=\Ad(q^{-2\theta(\la)}),\qquad
\theta(\la)=-\la+\rho-\frac{1}{2} \sum_{l}h_lh^l,
\]
where $\Ad(x)y=xyx^{-1}$. We def\/ine the face type twistor
$F(\la)\in U_q\widehat{\otimes} U_q$ as follows.

\begin{dfn}[Face type twistor]\lb{dfn:face}
\begin{gather}
F(\la)=\cdots \bigl((\varphi_\la)^2\otimes\id\bigr)\cR_0^{-1}
\bigl(\varphi_\la\otimes\id\bigr)\cR_0^{-1}
\nn\\
\phantom{F(\la)}= \mathop{\prod_{k\geq 1}}^{\kuru}
\Bigl((\varphi_\la)^k\otimes\id\Bigr)\cR_0^{-1},
\label{facetwistor}
\end{gather}
where $\mathop{\prod\limits_{k\geq 1}}\limits^{\kuru}A_k=\cdots
A_3A_2A_1$.
\end{dfn}

\noindent Note that the $k$-th factor in the product
\eqref{facetwistor} is a formal power series in
$x_i^k=q^{2k\langle \al_i,\la\rangle}$ $(i\in I)$ with leading
term $1$.

\begin{thm}[\cite{JKOS}]\lb{thm:face}
The twistor $F(\la)$ satisfies the shifted cocycle condition and
the norma\-li\-zation condition given by
\begin{gather}
1)\ \ F^{(12)}(\la)(\Delta\otimes \id)F(\la)
=F^{(23)}(\la+h^{(1)})(\id \otimes \Delta)F(\la),
\lb{facecocy}\\
2)\ \ \left(\epsilon\otimes\id\right)F(\la)
=\left(\id\otimes\epsilon\right)F(\la)=1. \lb{epF}
\end{gather}
\end{thm}

In \eqref{facecocy}, $\la$ and $h^{(1)}$ means
$\la=\sum_l\la_lh^l$
 and $h^{(1)}=\sum_l h_l^{(1)}h^l,\  h_l^{(1)}=h_l\otimes 1 \otimes 1$, respectively.
Note that from \eqref{strR}, one has
\[
[h\otimes 1+1\otimes h, F(\la)]=0\qquad \forall \, h\in \h.
\]

Now let us  def\/ine $\Delta_\la$, $\cR(\la)$, $\Phi(\la)$ and
$\al_\la$, $\beta_\la$ by
\begin{gather}
\Delta_{\la}(a)=F^{(12)}(\la) \,\Delta(a)\,  F^{(12)}(\la)^{-1},
\lb{facecopro}\\
\cR(\la)=F^{(21)}(\la)\,\cR\, F^{(12)}(\la)^{-1},
\lb{faceR}\\
\Phi(\la)=F^{(23)}(\la)F^{(23)}(\la+h^{(1)})^{-1},\\
\lb{facephi} \alpha_\la=\sum_i S(d_i)l_i, \qquad \beta_\la=\sum_i
m_i S(g_i)
\end{gather}
for $\sum_i k_i\otimes l_i=F(\la)^{-1}$, $\sum_i m_i \otimes
n_i=F(\la)$.

\begin{dfn}[Face type elliptic quantum group]\lb{dfn:Bqla}
With $S$ and $\epsilon$ defined by \eqref{coalgstr}, the set
 $(U_q(\goth{g}),\Delta_\la,S,\ve, \al_\la,\beta_\la, \Phi(\la),\cR(\la))$ forms
 a quasi-Hopf algebra~{\rm \cite{JKOS}}. We call it
the face type elliptic quantum group  $\Bqla(\goth{g})$.
\end{dfn}

From \eqref{YBE},  \eqref{facecocy} and \eqref{faceR}, one can
show that $\cR(\la)$ satisf\/ies the dynamical YBE.
\begin{thm}[Dynamical Yang--Baxter equation]
\begin{gather}
\cR^{(12)}(\la+h^{(3)})\cR^{(13)}(\la)\cR^{(23)}(\la+h^{(1)}) =
\cR^{(23)}(\la)\cR^{(13)}(\la+h^{(2)})\cR^{(12)}(\la). \lb{DYBE}
\end{gather}
\end{thm}

We hence call $\cR(\la)$ the universal dynamical $R$ matrix.

Now let us parametrize  $\la$ in the following way.
\begin{gather}
\la=\left(r+\frac{h^\vee}{ a_0^\vee}\right)d+s c+\bar{\la}\qquad
(r,s\in \C),\lb{parametlamda}
\end{gather}
where $\bar{\la}$ stands for the classical part of $\la$, and
$h^{\vee}$ denotes the dual Coxeter number of $\g$. Note also
$\rho=h^\vee \Lambda_0+\bar{\rho}$ and $d=a_0^\vee\Lambda_0$. Let
$\{\bar{h}_j\}$ and $\{\bar{h}^j(=\bar{\Lambda}_j)\}$ denote the
classical part of the  basis  and its dual of $\goth{h}$. We then
have
\begin{gather}
\varphi_\la=\Ad(p^dq^{2cd}q^{-2\bar{\theta}(\la)}), \qquad
\bar{\theta}(\la)=-\bar{\la}+\bar{\rho}-\frac{1}{2}\sum
\bar{h}_j\bar{h}^j.\lb{bartheta}
\end{gather}
Here we set $p=q^{2r}$. Set further
\begin{gather}
\cR(z)=\Ad(z^d\otimes 1)(\cR), \lb{Rg}
\\
F(z,\la)=\Ad(z^d\otimes 1)(F(\la)), \lb{Fg}
\\
\cR(z,\la)=\Ad(z^d\otimes 1)(\cR(\la))
=\sigma(F(z^{-1},\la))\cR(z)F(z,\la)^{-1}. \lb{Rg2}
\end{gather}
Then $\cR(z)$ and $F(z,\la)$ are formal power series in $z$,
whereas $\cR(z,\la)$ contains both positive and negative powers of
$z$.

 From the def\/inition \eqref{facetwistor} of $F(\la)$, one can easily derive
 the following dif\/ference equation for the twistor.
\begin{thm}[Dif\/ference equation \cite{JKOS}]\lb{DeqF}
\begin{gather}
F(pq^{2c^{(1)}}z,\la)=\Ad(q^{2\bar{\theta}(\la)}\otimes\id)
\bigl(F(z,\la)\bigr)\cdot q^T\cR(pq^{2c^{(1)}}z). \lb{Fdiff}
\end{gather}
\end{thm}

Furthermore, noting $\Ad(z^d)(e_i)=z^{\delta_{i,0}}e_i$, one can
drop all the $e_0$ dependent terms in $\cR(z)$ and $F(z,\la)$ by
taking  the limit $z\to 0$. We thus obtain
\begin{gather}
\lim_{z\to 0}q^{c\otimes d+d\otimes c}\cR(z)=\cR_{\bar{\goth{g}}},\\
\lim_{z\to 0}F(z,\la)=F_{\bar{\goth{g}}}(\bar{\la}),\lb{Fini}
\end{gather}
where $\cR_{\bar{\goth{g}}}$ and $F_{\bar{\goth{g}}}(\bar{\la})$
are the universal $R$ matrix and the twistor of
$U_q(\bar{\goth{g}})$. Then from \eqref{Fdiff}, we obtain the
following equation for $F_{\bar{\goth{g}}}(\bar{\la})$.
\begin{lem}
\begin{gather}
F_{\bar{\goth{g}}}(\la)=\Ad(q^{2\bar{\theta}(\la)}\otimes\id)
\bigl(F_{\bar{\goth{g}}}(\bar{\la})\bigr)\cdot
q^{\bar{T}}\cR_{\bar{\goth{g}}},\lb{ABRReq}
\end{gather}
where $\bar{T}=\sum\limits_{i=1}^n \bar{h}_i\otimes \bar{h}^i$.
\end{lem}

\noindent {\bf Remark.} \eqref{ABRReq} corresponds to (18) in
\cite{ABRR},
 where the comultiplication and the universal $R$ matrix
 are our $\Delta^{op}$ and $\cR^{-1}_{\bar{\goth{g}}}$, respectively.

\begin{lem}[\cite{ABRR}]\lb{ABRR}
The equation \eqref{ABRReq} has the unique solution
$F_{\bar{\goth{g}}}(\bar{\la})\in
U_q(\bar{\gb}_+)\widehat{\otimes} U_{q}(\bar{\gb}_-)$ in the form
$F_{\bar{\goth{g}}}(\bar{\la})=1+\cdots$. Here $U_q(\bar{\gb}_+)$
(resp. $U_{q}(\bar{\gb}_-)$) is the subalgebra of $U_q(\bar{\g})$
generated by $e_i$, $t_i$ $(i=1,2,\dots,n)$
 (resp. $f_i$, $t_i$ $(i=1,2,\dots,n)$).
\end{lem}

\begin{thm}
For $\la\in \h$ given by \eqref{parametlamda}, the difference
equation \eqref{Fdiff} has a unique solution.
\end{thm}

\begin{proof}
Let us set $\bar{\varphi}_\la=\Ad(q^{-2\bar{\theta}(\la)})$.
Iterating \eqref{Fdiff}, $N$ times, we obtain
\[
F(z,\la)=\bigl(\bar{\varphi}^N_\la\otimes\id\bigr)F((pq^{2c^{(1)}})^Nz,{\la})
\mathop{\prod_{N\geq k\geq 1}}^{\kuru}
\Bigl((\bar{\varphi}_\la)^k\otimes\id\Bigr)\cR_0((pq^{2c^{(1)}})^kz)^{-1}.
\]
Taking the limit $N \to \infty$, one obtains
\[
F(z,\la)= A \mathop{\prod_{ k\geq 1}}^{\kuru}
\bigl((\bar{\varphi}_\la)^k\otimes\id\bigr)\cR_0((pq^{2c^{(1)}})^kz)^{-1},
\]
where we set
\begin{gather*}
A=\lim_{N\to \infty}
\bigl(\bar{\varphi}^N_\la\otimes\id\bigr)F((pq^{2c^{(1)}})^Nz,{\la})\\
\phantom{A}{}=\lim_{N\to \infty}
\bigl(\bar{\varphi}^N_\la\otimes\id\bigr){F}_{\bar{\goth{g}}}(\bar{\la}).
\end{gather*}
Then the statement follows from Lemma~\ref{ABRR}.
\end{proof}

\section{Vertex operators and  fusion matrices}

\subsection[The vertex operators of $U_q(\gbig)$]{The vertex operators of $\boldsymbol{U_q(\gbig)}$}

Let $V$ and $W$ be f\/inite dimensional irreducible modules of
$U'_q$. Let $\la, \mu \in \h^*$ be level-$k$ generic  elements
such that  $\langle c,\la\rangle=\langle c,\mu\rangle =k$. We
denote by $M_\la$ and $M_\mu$ the two irreducible Verma modules
with highest weights $\la$ and $\mu$, respectively.

\begin{df}[Vertex operator]
Writing
$\triangle_\la=\frac{(\la|\la+2\rho)}{2(k+h^\vee)}$\footnote{Hopefully,
there is no confusion of $\triangle_\la$ with $\Delta_\la$ in
\eqref{facecopro}.}, let us consider the formal series  given by
\begin{gather}
\Psi^{\mu}_{\la}(z)=z^{\triangle_\mu-\triangle_\la}\tPsi^{\mu}_{\la}(z),\qquad
\tPsi^{\mu}_{\la}(z)=\sum_{j}\sum_{n\in \Z}v_jz^{-n}\otimes
(\tPsi^{\mu}_{\la})_{j,n}. \lb{typeIIVO}
\end{gather}
Here $\{v_j\}$ denotes a weight basis of $V$. The coefficients $
(\tPsi^{\mu}_{\la})_{j,n}$ are the maps
\begin{gather}
(\tPsi^{\mu}_{\la})_{j,n}: (M_\la)_\xi \to
(M_\mu)_{\xi-\wt(v_j)+n\delta}, \lb{Psiweight}
\end{gather}
such that $\tPsi^{\mu}_{\la}(z)$ is the $U_q$-module intertwiners
\begin{gather}
\tPsi^{\mu}_{\la}(z)\ :\ M_\la \to V_z \widehat{\otimes}\ M_\mu,\nn \\
\tPsi^{\mu}_{\la}(z)\ x=\Delta(x)\ \tPsi^{\mu}_{\la}(z)\qquad
\forall \, x\in U_q.\lb{Psiint}
\end{gather}
Here $\widehat{\otimes}$ denotes a formal completion
\[
M\widehat{\otimes}\ N=\bigoplus_{\nu} \prod_{\xi} M_\xi {\otimes}\
N_{\xi-\nu}.
\]
We call $\Psi^{\mu}_{\la}(z)$ the vertex operator (VO) of $U_q$.
\end{df}

\noindent {\bf Remark.} $\Psi^{\mu}_{\la}(z)$ is the type II VO in
the terminology of~\cite{JM}.

We also def\/ine  $U'_q$-module intertwiners $\Psi^{\mu}_{\la} :
M_\la\to V\otimes \widehat{M}_\mu$ by
\begin{gather}
\Psi^{\mu}_{\la}=\sum_{j}\ v_j\otimes \left(\sum_{n\in \Z}
(\tPsi^{\mu}_{\la})_{j,n}\right).\lb{ptypeIIVO}
\end{gather}
Here $\widehat{M}_\mu=\prod_{\nu}(M_\mu)_\nu$. Note that there is
a bijective correspondence between $\Psi^{\mu}_{\la}(z)$ and
$\Psi^{\mu}_{\la}$.

Let $u_\la$ and $u_\mu$ denote the highest weight vectors of
$M_\la$ and $M_\mu$, respectively. Let us write the image of
$u_\la$ by the VO as
\begin{gather}
\Psi^{\mu}_{\la}\ u_{\la}=v\otimes u_\mu+\sum_{i'}v_{i'}\otimes
u_{i'}, \lb{leadingtermII}
\end{gather}
where $u_{i'}\in M_\mu$, $\wt(u_{i'})< \mu$ and $v, v_{i'}\in V$.
We call the vector $v$ the leading term of $\Psi^\mu_\la$. Note
that from \eqref{Psiweight},
$cl(\la)=\wt(v)+cl(\mu)=\wt(v_{i'})+cl(\wt(u_{i'}))$. We set
\[
V^\mu_\la=\{ v\in V\ |\ \wt(v)=cl(\la-\mu) \}.
\]

\begin{thm}[\cite{FR,DJO}]\lb{DJO}
The map $ \langle \  \rangle : \Psi^\mu_\la \mapsto \langle
\id\otimes u_\mu^*,\Psi^\mu_\la u_\la\rangle$ gives a
$\C(q^{1/2})$-linear isomorphism
\[
\Hom_{U'_q(\g)}(M_\la,V \otimes
\widehat{M}_\mu)\stackrel{\sim}{\longrightarrow} V^\mu_\la
\]
\end{thm}

This theorem tells that $\Psi^\mu_\la$ is determined by its
leading term. Namely, for given $v_{0} \in V^\mu_\la$, there
exists a unique VO satisfying
\[
\langle \Psi^\mu_\la\rangle=v_{0}.
\]
We denote such VO by $\Psi^{\mu, v_{0}}_\la$ and corresponding
$U_q$-intertwiner by $\Psi^{\mu, v_{0}}_\la(z)$.

\begin{prop}[\cite{DJO}]
Let $\{v_j\}$ be a basis of $V^\mu_\la$. The set of VOs
$\{\Psi^{\mu, v_j}_\la\}$ forms a basis of $\Hom_{U'_q(\g)}(M_\la,
V\otimes\widehat{M}_\mu )$.
\end{prop}

\subsection[The $q$-KZ equation and connection matrices]{The $\boldsymbol{q}$-KZ equation and connection matrices}

Let $\la, \mu, \nu\in \h^*$ be  level-$k$ elements. Let $\{v_i\}$
and  $\{w_j\}$ be weight bases of $V_\la^\mu$ and $W_\mu^\nu$,
respectively. Consider the VOs $\Psi^{\mu, v_{i}}_\la$ and
$\Psi^{\nu, w_{j}}_\mu$ given by
\begin{gather*}
\Psi^{\mu, w_{i}}_\la(z_1)\ :\ M_\la \to W_{z_1}\widehat{\otimes}\
M_\mu
,\\
\Psi^{\nu, v_{j}}_\mu(z_2) \ :\ M_\mu \to
V_{z_2}\widehat{\otimes}\ M_\nu,
\end{gather*}
and their composition
\[
\left(\id\otimes \Psi^{\nu, v_{i}}_\mu(z_2)\right)\Psi^{\mu,
w_{j}}_\la(z_1)\ :\ M_\la \to W_{z_1}{\otimes}\  V_{z_2}
\widehat{\otimes}\ {M_\nu}.
\]
Setting
\[
\Psi^{(\nu,\mu,\la)}(z_1,z_2)=\left< \id\otimes \id\otimes
u_\nu^*, \left(\id\otimes \Psi^{\nu,
v_{i}}_\mu(z_2)\right)\Psi^{\mu, w_{j}}_\la(z_1)u_\la\right>,
\]
we call $\Psi^{(\nu,\mu,\la)}(z_1,z_2)$ the two-point function.

\begin{thm}[$\boldsymbol{q}$-KZ equation \cite{FR,IIJMNT}]\lb{qKZ}
The two-point function $\Psi^{(\nu,\mu,\la)}(z_1,z_2)$ satisfies
 the $q$-KZ equation
\begin{gather}
\Psi^{(\nu,\mu,\la)}(q^{2(k+h^{\vee})}z_1,z_2)=(q^{-\pi_W(\bar{\nu}+\bar{\la}+2\bar{\rho})}
\otimes \id)R_{WV}(z)\Psi^{(\nu,\mu,\la)}(z_1,z_2),\lb{qKZeq}
\end{gather}
where $R_{WV}(z)=(\pi_W\otimes \pi_V)\cR(z)$.
\end{thm}

\begin{proof}
 See Appendix.
 \end{proof}

\noindent {\bf  Remark~\cite{FR}.} A solution of the $q$-KZ
equation \eqref{qKZeq} is a function of  $z=z_1/z_2$ and has
a~form $G(z)f(z_1,z_2)$. Here $G(z)$ is a meromorphic function
multiplied by a fractional power of $z$ determined from the
normalization function of the $R$ matrix $R_{WV}(z)$, while
$f(z_1,z_2)$ is an analytic function in $|z_1|>|z_2|$ and can be
continued meromorphically to
 $(\C^{\times})^2$. Hence a~solution of the $q$-KZ equation is uniquely determined,
if one f\/ixes the normalization. It also follows that the
composition of the VOs is well def\/ined in the region
$|z_1|>|z_2|$ and  can be continued meromorphically to
 $(\C^{\times})^2$ apart from an overall fractional power of $z$.

The following commutation relation holds in the sense of analytic
continuation.

\begin{thm}[Connection formula \cite{FR,DJO}]\lb{CommRel}
\begin{gather*}
\left( PR_{WV}(z_1/z_2)\otimes \id \right)\left(\id\otimes \Psi^{\nu v_{i}}_\mu(z_2)\right)\Psi^{\mu w_{j}}_\la(z_1)\\
\qquad{}=\sum_{i', j', \mu'}\left(\id\otimes \Psi^{\nu
w_{j'}}_{\mu'}(z_1)\right) \Psi^{\mu' v_{i'}}_\la(z_2)
C_{WV}\left( \left. \mmatrix{\la&w_{j}&\mu\cr
           v_{i'}& &v_{i}\cr
           \mu'&w_{j'}&\nu\cr} \right|\frac{z_1}{z_2}\right),
\end{gather*}
where $v_{i'}$ and $w_{j'}$ are base vectors of $V_\la^{\mu'}$ and
$W_{\mu'}^\nu$, respectively.
\end{thm}

The matrix $C_{WV}$ is called the connection matrix. The following
theorem states basic pro\-per\-ties of the connection matrix.

\begin{thm}[\cite{FR,DJO}]\lb{connectMat}
\begin{itemize}\itemsep=0pt
\item[{\rm 1)}] The matrix elements of $C_{WV}$ are given by a
ratio of elliptic theta functions.

\item[{\rm 2)}] The matrix $C_{WV}$ satisfies i) the face type YBE
and ii) the unitarity condition (the first inversion relation):
\begin{gather*}
i)\quad \sum_{v_l,w_i,u_j,\mu''} C_{VU}\left( \left.
\mmatrix{\la&u_{j}&\mu''\cr
           w_{i'}& &w_{i}\cr
           \omega&u_{j'}&\nu\cr} \right|\frac{z_2}{z_3}\right)
C_{WU}\left( \left. \mmatrix{\la'&u_{j''}&\mu'\cr
           v_{l'}& &v_{l}\cr
           \la&u_{j}&\mu''\cr} \right|\frac{z_1}{z_3}\right)
\\
\qquad \quad {}\times C_{WV}\left( \left.
\mmatrix{\mu'&w_{i''}&\mu\cr
           v_{l}& &v_{l''}\cr
           \mu''&w_{i}&\nu\cr} \right|\frac{z_1}{z_2}\right)\\
\qquad{}=\sum_{v_l,w_i,u_j,\mu''} C_{WV}\left( \left.
\mmatrix{\la'&w_{i}&\mu''\cr
           v_{l'}& &v_{l}\cr
           \la&w_{i'}&\omega\cr} \right|\frac{z_1}{z_2}\right)
C_{WU}\left( \left. \mmatrix{\mu''&u_{j}&\mu\cr
           v_{l}& &v_{l''}\cr
           \omega&u_{j'}&\nu\cr} \right|\frac{z_1}{z_3}\right)
\\
\qquad \quad{}\times C_{VU}\left( \left.
\mmatrix{\la'&u_{j''}&\mu'\cr
           w_{i}& &w_{i''}\cr
           \mu''&u_{j}&\mu\cr} \right|\frac{z_2}{z_3}\right),\\
ii)\quad\sum_{v_{i'}, w_{j'}, \mu'} C_{WV}\left( \left.
\mmatrix{\la&w_{j}&\mu\cr
           v_{i'}& &v_{i}\cr
           \mu'&w_{j'}&\nu\cr} \right|z\right)
C_{VW}\left( \left. \mmatrix{\la&v_{i'}&\mu'\cr
           w_{j''}& &w_{j'}\cr
           \mu''&v_{i''}&\nu\cr} \right|z^{-1}\right)\\
\qquad{}=\delta_{w_j,w_{j''}}\delta_{v_i,v_{i''}}\delta_{\mu,\mu''}.
\end{gather*}
\item[{\rm 3)}] In the case $V=W$, $C_{VV}$ satisfies the second
inversion relation
\begin{gather*}
\sum_{v_{i'},w_{j'},\la} \frac{G_\la G_\nu}{G_\mu G_{\mu'}}
C_{VV}\left( \left. \mmatrix{\la&w_{j}&\mu\cr
           v_{i'}& &v_{i}\cr
           \mu'&w_{j'}&\nu\cr} \right|z^{-1}\right)
C_{VV}\left( \left. \mmatrix{\la&v_{i'}&\mu'\cr
           w_{j}& &w_{j''}\cr
           \mu&v_{i''}&\nu'\cr} \right|\xi^{-2}z\right)\\
           \qquad{} =\alpha_{VV}(z)\delta_{v_i,v_{i''}}\delta_{w_{j},w_{j''}}\delta_{\nu,\nu'}.
\end{gather*}
Here $\xi=q^{t h^\vee}$ ($t=({\rm long\ root })^2/2$), and
$G_\la$ and $\alpha_{VV}(z)$ are given by {\rm (5.16)} and {\rm
(5.13)} in~{\rm \cite{DJO}}, respectively.
\end{itemize}
\end{thm}

Furthermore, if $V_z$ is self dual i.e.\ there exists
 an isomorphism of $U_q$-modules $Q: V_{\xi^{-1} z}\simeq V_z^{*a}$,
we have the following crossing symmetry.
\begin{thm}[\cite{DJO}]
\begin{gather}
\beta_{VV}(z^{-1})C_{VV}\left( \left. \mmatrix{\mu&w_{j}&\nu\cr
           v_{i}& &w_{j'}\cr
           \la&v_{i'}&\mu'\cr} \right|\xi^{-1}z^{-1}\right)
=\sum_{\tilde{i},\tilde{j}}\gamma^{\mu i\tilde{i}}_\la g^{\nu
\tilde{j}j'}_{\mu'} C_{VV}\left( \left.
\mmatrix{\la&v_{\tilde{i}}&\mu\cr
           v_{i'}& &w_{j}\cr
           \mu'&w_{\tilde{j}}&\nu\cr} \right|z\right), \lb{cross}
\end{gather}
where $g^{\nu jj'}_{\mu}$ denotes a certain matrix element
appearing in the inversion relation of the VO's, and $\gamma^{\nu
}_{\mu}$ is its inverse matrix such that $\sum_{\tilde{j}}g^{\nu
j\tilde{j}}_{\mu}\gamma^{\nu \tilde{j}j'}_{\mu}=\delta_{jj'}$.
\end{thm}

In Section~5, we will discuss the evaluation molude $V_z$ with $V$
being the vector representation for $\g=A_n^{(1)}, B_n^{(1)},
C_n^{(1)}, D_n^{(1)}$. There $V_z$ is self dual except for
$A_n^{(1)}$ $(n>1)$, and  is multiplicity free. Therefore, $\dim
V^\mu_\la=1$ etc. Hence  the matrices $g^\nu_\mu, \gamma^{\nu
}_{\mu}$ are scalars satisfying $\gamma^{\nu }_{\mu}=1/g^\nu_\mu$.
In this case, let us consider the gauge transformation
\[
C_{VV}\left( \left. \mmatrix{\mu&\nu\cr
           \la&\mu'\cr} \right|z\right)=f(z)
\frac{F(\mu,\nu)F(\nu,\mu')}{F(\mu,\la)F(\la,\mu')}\tilde{C}_{VV}\left(
\left. \mmatrix{\mu&\nu\cr
           \la&\mu'\cr} \right|z\right)
\]
with the choice
\begin{gather*}
f(z)f(z^{-1})=1, \qquad f(\xi^{-1}z^{-1})=\beta_{VV}(\xi z)f(z),\\
F(\nu,\mu)F(\mu,\nu)=g^\nu_\mu\sqrt{\frac{G_\mu}{G_\nu}}.
\end{gather*}
Then we can change the crossing symmetry relation \eqref{cross} to
\begin{gather*}
\tilde{C}_{VV}\left( \left. \mmatrix{\mu&\nu\cr
           \la&\mu'\cr} \right|\xi^{-1}z^{-1}\right)
=\sqrt{\frac{G_\nu G_\la}{G_\mu G_{\mu'}}} \tilde{C}_{VV}\left(
\left. \mmatrix{\la&\mu\cr
           \mu'&\nu\cr} \right|z\right). \lb{Gcross}
\end{gather*}
The same gauge transformation makes  the face type YBE $i)$, the
unitarity condition $ii)$ and the second inversion relation $3)$
in Theorem \ref{connectMat} unchanged except for the factor
$\alpha_{VV}(z)$ in the RHS of $3)$, which is changed to 1.

\subsection{Fusion matrices}

We here follows the idea by Etingof and Varchenko \cite{EV2},
where the cases $U_q(\bar{\g})$ with $\bar{\g}$ being simple Lie
algebras are discussed. We extend their results to the cases
$U_q(\g)$ with $\g$ being af\/f\/ine Lie algebras.

\begin{df}[Fusion matrix]
Fix $\la\in \h^*$. The  fusion matrix is defined to be a
$\h$-linear map $ J_{WV}(\la) : W\otimes V \to W\otimes V$
satisfying
\begin{gather*}
J_{VW}(\la)=\bigoplus_{\nu} J_{VW}(\la)_\nu,\\
J_{WV}(\la)_\nu(w_{j}\otimes v_{i})= \left<\id\otimes \id\otimes
u_\nu^*, \left(\id\otimes \Psi^{\nu, v_{i}}_\mu\right)\Psi^{\mu,
w_{j}}_\la u_\la\right>\ \in
 (W\otimes V)_{cl(\la-\nu)},
\end{gather*}
for $v_{i}\in V^\nu_\mu$, $w_{j}\in W^\mu_\la$.
\end{df}

Note that from \eqref{Psiweight},
\begin{gather}
[h\otimes 1+1\otimes h, J_{WV}(\la)]=0 \qquad \forall\, h\in
\h.\lb{hlinear}
\end{gather}

Noting \eqref{leadingtermII} and the intertwining property of the
vertex operators, we have
\begin{gather}
J_{WV}(\la)_\nu(w_{j}\otimes v_{i}) = w_{j}\otimes
v_{i}+\sum_{l}C_l(\la)w_l\otimes v_l, \lb{fusionmat}
\end{gather}
where $\wt(v_{l})< \wt(v_{i})$, and $C_l(\la)$ is a function of
$\la$. Hence $J_{VW}(\la)$ is an upper triangular matrix with all
the diagonal components being 1. Therefore
\begin{thm}
The fusion matrix $J_{WV}(\la)$ is invertible.
\end{thm}

The def\/inition of $J_{WV}(\la)$ indicates that the leading term
of the intertwiner $\left(\id\otimes \Psi^{\nu,
v_{i}}_\mu\right)\Psi^{\mu, w_{j}}_\la$ is $J_{WV}(\la)(
w_{j}\otimes v_{i})$. Therefore we write
\[
\Psi^{\nu, J_{WV}(\la)(w_{j}\otimes v_i)}_\la=\left(\id\otimes
\Psi^{\nu, v_{i}}_\mu\right)\Psi^{\mu, w_{j}}_\la.
\]
Let us def\/ine $J_{U, W\otimes V}(\omega)$ and  $J_{U\otimes W,
V}(\omega):U\otimes W\otimes V\to U\otimes W\otimes V$ for $\omega
\in \h^*$ by
\begin{gather*}
J_{U, W\otimes V}(\omega)(u_l\otimes w_j\otimes v_i)= \big\langle
\id\otimes \id\otimes \id\otimes u^*_{\nu}, \big(\id\otimes
\Psi^{\nu, J_{WV}(\la)(w_j\otimes v_i)}_{\la}
\big)\Psi^{\la, u_{l}}_\omega u_\omega \big\rangle, \\
J_{ U\otimes W, V}(\omega)(u_l\otimes w_j\otimes
v_i)=\big\langle\id\otimes \id\otimes \id\otimes u^*_{\nu} ,
(\id\otimes \id\otimes \Psi^{\nu, v_i}_{\mu})\Psi^{\mu, J_{WU}
 (\omega)(w_{j} \otimes u_l)}_\omega u_\omega\big\rangle,
\end{gather*}
where $u_l\in U_\omega^\la$.

\begin{thm}
The fusion matrix satisfies the shifted cocycle condition
\begin{gather}
J_{U\otimes W, V}(\omega)(J_{UW}(\omega)\otimes \id)= J_{U,
W\otimes V}(\omega)(\id\otimes J_{WV}(\omega-h^{(1)}) )\qquad {\rm
on }\quad U\otimes W\otimes V, \lb{Jshiftcocycle}
\end{gather}
where $h v_j=\wt(v_j) v_j \ (v_j\in V)$ etc.
\end{thm}

\begin{proof} Consider the composition
\begin{gather*}
(\id\otimes \id\otimes\Psi^{\nu, v_i}_\mu)(\id\otimes \Psi^{\mu, w_j}_\la)\Psi^{\la, u_l}_\omega:\\
M_\omega \stackrel{\Psi^{\la, u_l}_\omega}{\longrightarrow}
U\otimes  \widehat{M}_\la \stackrel{\id\otimes \Psi^{\mu,
w_j}_\la}{\longrightarrow}  U\otimes W  \otimes\widehat{M}_\mu
\stackrel{ \id\otimes \id\otimes\Psi^{\nu,
v_i}_\mu}{\longrightarrow}  U\otimes W \otimes V
\otimes\widehat{M}_\nu
\end{gather*}
and express the highest to highest expectation value of it in two
ways, and use $cl(\la)=cl(\omega)-\wt(u_l)$.
\end{proof}

\noindent {\bf Remark.} Regarding $(\Delta\otimes
\id)J(\omega)=J_{U\otimes W, V}(\omega)$, $(\id\otimes
\Delta)J(\omega)=J_{U, W\otimes V}(\omega)$,
$J^{(12)}(\omega)=J_{UW}(\omega)\otimes \id$ and
$J^{(23)}(\omega)=\id\otimes J_{WV}(\omega)$, one obtains
\eqref{facecocy} from \eqref{Jshiftcocycle} by replacing
$J(\omega)$ by $F^{-1}(-\omega)$.

Now let $\la\in \h^*$  be a level-$k$ element.
 By using the $U_q$-module VOs \eqref{typeIIVO},
we def\/ine a $\h$-linear map
 $J_{WV}(z_1,z_2;\la)\  :\ W_{z_1}\otimes V_{z_2} \to W_{z_1}\otimes V_{z_2}$
by
\begin{gather*}
J_{WV}(z_1,z_2;\la)=\bigoplus_{\nu} J_{WV}(z_1,z_2;\la)_\nu,\\
J_{WV}(z_1,z_2;\la)_\nu(w_{j}\otimes v_{i})\\
\qquad{}= \big\langle \id\otimes \id\otimes u_\nu^*,
\big(\id\otimes \tPsi^{\nu, v_{i}}_\mu(z_2)\big)\tPsi^{\mu,
w_{j}}_\la(z_1) u_\la\big\rangle \in
 (W_{z_1}\otimes V_{z_2})_{cl(\la-\nu)},
\end{gather*}
for $v_{j}\in V^\nu_\mu$, $w_{j}\in W^\mu_\la$. Then from the
$q$-KZ equation \eqref{qKZeq}, one can derive  the following
dif\/ference equation for $J_{WV}(z_1,z_2;\la)$.

\begin{lem}\lb{DeqK}
\begin{gather}
J_{WV}(q^{2(k+h^\vee)}z_1,z_2; \la)(q^{-2\pi_W(\bar{\theta}(-\la))}\otimes \id)\nn\\
\qquad{}=(q^{-2\pi_W(\bar{\theta}(-\la))} \otimes \id )
q^{\pi_{W\otimes V}(\bar{T})}R_{WV}(z_1/z_2)  J_{WV}(z_1,z_2,\la).
\lb{DiffK}
\end{gather}
\end{lem}

\begin{proof}
 See Appendix.
 \end{proof}

From the remark below Theorem \ref{qKZ}, $J_{WV}(z_1,z_2; \la)$ is
a function of the ratio $z=z_1/z_2$. Let us parameterize a
level-$k$ $\la\in\h^*$ as \eqref{parametlamda}. Comparing Theorem
\ref{DeqF} and Lemma \ref{DeqK}, we f\/ind that the dif\/ference
equation for $F_{WV}(z,-\la)^{-1}=(\pi_W\otimes
\pi_V)F(z,-\la)^{-1}$ coincides with  the $q$-KZ equation
\eqref{DiffK} for $J_{WV}(z_1,z_2;\la)$ on $w_j\otimes v_i$ under
the identif\/ication $r=-(k+h^\vee)$. Hence the uniqueness of the
solution to the $q$-KZ equation yields
 the following theorem.
\begin{thm}\lb{PsiPsiF}
For a level-$k$  $\la\in \h^*$ in the parameterization
\eqref{parametlamda},
\[
\big\langle \id\otimes \id\otimes u_{\nu}^*,\big(\id \otimes
\tPsi^{\nu,v_i}_\mu(z_2)\big)
\tPsi^{\mu,w_j}_{\la}(z_1)u_{\la}\big\rangle
=F_{WV}\left({z_1}/{z_2},-\la\right)^{-1}(w_j\otimes v_i).
\]
\end{thm}

\noindent {\bf Remark.} Relation between the twistors and the
fusion matrices was f\/irst discussed by Etingof and Varchenko for
the case $\g$ being simple Lie algebra (Appendix~9 in~\cite{EV2}).
 Their coproduct and the universal $R$ matrix correspond
 to  our $\Delta^{op}$ and $\cR^{-1}_{\bar{\g}}$, respectively,
 and the twistor is identif\/ied with the two point
 function of the $U_q(\bar{\g})$-analogue of the type I VOs.

\section[Dynamical $R$ matrices and connection matrices]{Dynamical $\boldsymbol{R}$ matrices and connection matrices}

Let $(\pi_V,V)$, $(\pi_W, W)$ be f\/inite dimensional
representations of $U'_{q}(\g)$, and
 $\{v_i\}$, $\{w_j\}$ be their weight bases, respectively.
Now we consider the dynamical $R$ matrices given as the images of
the universal $R$ matrix $\cR(\la)$
\begin{gather}
R_{WV}(z,\la)=(\pi_W\otimes \pi_V)\cR(z,\la)\nn\\
\phantom{R_{WV}(z,\la)}{}=F_{VW}^{(21)}(z^{-1},\la)R_{WV}(z)F_{WV}(z,\la)^{-1}.\lb{DynamicalR}
\end{gather}
Note that  $F^{(21)}_{VW}(z^{-1},\la)=PF_{VW}(z^{-1},\la)P$.

By using Theorems \ref{PsiPsiF} and \ref{CommRel}, we show that
the dynamical $R$ matrix $R_{WV}(z,\la)$ coincides with
 the corresponding connection matrix for the $q$-KZ equation  of $U_q(\g)$ associated with the representations $(\pi_V,V)$ and $(\pi_W, W)$.
\begin{thm}\lb{RmatCmat} For level-$k$
 $\la\in \h^*$ in the form \eqref{parametlamda} with $r=-(k+h^\vee)$,  we have
\[
R_{WV}(z,-\la)(w_j\otimes v_i)= \sum_{i',j',\mu'} C_{WV}\left(
\left. \mmatrix{\la&w_{j}&\mu\cr
           v_{i'}& &v_{i}\cr
           \mu'&w_{j'}&\nu\cr} \right|z\right)(w_{j'}\otimes v_{i'}),
\]
where $v_i\in V_\mu^\nu$, $w_j\in W_\la^\mu$, $v_{i'}\in
V_{\la}^{\mu'}$ and $w_{j'}\in W_{\mu'}^\nu$.
\end{thm}

\begin{proof}
 From Theorem \ref{PsiPsiF} and \eqref{DynamicalR}, we have
\begin{gather*}
R_{WV}(z_1/z_2,-\la)(w_j\otimes v_i)\\
\qquad{}=F_{VW}^{(21)}(z^{-1},-\la)R_{WV}(z)\big\langle \id\otimes
\id\otimes u_{\nu}^*, \big(\id \otimes
\tPsi^{\nu,v_i}_\mu(z_2)\big)
\tPsi^{\mu,w_j}_{\la}(z_1)u_{\la}\big\rangle,
\end{gather*}
where we set $z=z_1/z_2$. Then from Theorem \ref{CommRel}, we
obtain
\begin{gather*}
R_{WV}(z_1/z_2,-\la)(w_j\otimes v_i)=F_{VW}^{(21)}(z^{-1},-\la)\\
\qquad{}\times \sum_{i',j',\mu'}P\big\langle \id\otimes \id\otimes
u_{\nu}^*,\big(\id \otimes \tPsi^{\nu,w_{j'}}_{\mu'}(z_1)\big)
\tPsi^{\mu,v_{i'}}_{\la}(z_2)u_{\la}\big\rangle C_{WV}\left(
\left. \mmatrix{\la&w_{j}&\mu\cr
           v_{i'}& &v_{i}\cr
           \mu'&w_{j'}&\nu\cr} \right|z\right)\\
\qquad{}=F_{VW}^{(21)}(z^{-1},-\la)
\sum_{i',j',\mu'}PF(z_2/z_1,-\la)^{-1}(v_{i'}\otimes
w_{j'})C_{WV}\left( \left. \mmatrix{\la&w_{j}&\mu\cr
           v_{i'}& &v_{i}\cr
           \mu'&w_{j'}&\nu\cr} \right|z\right)\\
\qquad{}=\sum_{i',j',\mu'}C_{WV}\left( \left.
\mmatrix{\la&w_{j}&\mu\cr
           v_{i'}& &v_{i}\cr
           \mu'&w_{j'}&\nu\cr} \right|z\right)(w_{j'}\otimes v_{i'}).\tag*{\qed}
\end{gather*}
 \renewcommand{\qed}{}
\end{proof}

Note that in view of Theorem~\ref{connectMat}, this theorem
indicates that the dynamical $R$ matrices of $\Bqla(\g)$ are
elliptic.

\section{Vector representations}
In this section, we consider the vector representation of the
universal $R$ matrix $\cR(\la)$ of $\Bqla(\g)$ for  $\g=A^{(1)}_n,
B^{(1)}_n, C^{(1)}_n, D^{(1)}_n$, and show that they coincide with
the corresponding face weights obtained by Jimbo, Miwa and
Okado~\cite{JMO}.

\subsection[Jimbo-Miwa-Okado's solutions]{Jimbo--Miwa--Okado's solutions}

Let us summarize Jimbo--Miwa--Okado's elliptic solutions to the
face type YBE.

Let $(\pi_V,V)$ be the vector representation of $U'_q(\g)$. We set
$\dim V=N$. Then $N=n+1,\, 2n+1$, $2n,\, 2n$ for $\g=A^{(1)}_n,
B^{(1)}_n, C^{(1)}_n, D^{(1)}_n$, respectively. Let us def\/ine an
index set $J$ by
\begin{gather*}
J=\{1,2,\dots,n+1\}\qquad {\rm for }\ A^{(1)}_n\\
\phantom{J}{}=\{0,\pm1,\dots,\pm n\}\qquad {\rm for }\ B^{(1)}_n\\
\phantom{J}{}=\{\pm1,\dots,\pm n\}\qquad {\rm for }\ C^{(1)}_n,
D^{(1)}_n
\end{gather*}
and introduce a linear order $\prec$ in $J$ by
\[
1\prec 2\prec \cdots \prec n \ (\prec 0)\ \prec -n\prec \cdots
\prec -2\prec -1.
\]
We  also use the usual order $<$ in $J$.

Let $\bar{\Lambda}_j$ $(1\leq j\leq n)$ be the fundamental weights
of $\bar{\g}$. Following Bourbaki~\cite{Bourbaki} we introduce
orthonormal vectors $\{ \vep_1,\dots,\vep_n \} $ with the bilinear
form $(\vep_i|\vep_j)=\delta_{i,j}$. Then one has the following
expression of $\bLa_j$ as well as the set $\cA$ of weights
belonging to the vector representation of $\bar{\g}$.
\begin{alignat*}{3}
& A_n:\quad &&\cA=\{\vep_1-\vep, \ldots,\vep_{n+1}-\vep\},&\\
& &&\bLa_j=\vep_1+\cdots+\vep_j-j\vep\quad (1\leq j\leq n), \qquad
\vep=\frac{1}{n+1}
\sum_{j=1}^{n+1}\vep_j,&\\
& B_n:&&\cA=\{\pm\vep_1, \ldots,\pm\vep_{n},0\},&\\
&&&\bLa_j=\vep_1+\cdots+\vep_j\qquad (1\leq j\leq n-1),&\\
&&&\phantom{\bLa_j}{} =\frac{1}{2} (\vep_1+\cdots+\vep_n)\qquad (j=n),&\\
& C_n:&&\cA=\{\pm\vep_1, \ldots,\pm\vep_{n}\},&\\
&&&\bLa_j=\vep_1+\cdots+\vep_j\qquad (1\leq j\leq n),&\\
&D_n:&&\cA=\{\pm\vep_1, \ldots,\pm\vep_{n}\},&\\
&&&\bLa_j=\vep_1+\cdots+\vep_j\qquad (1\leq j\leq n-2),&\\
&&&\phantom{\bLa_j}{} =\frac{1}{2} (\vep_1+\cdots+\vep_{n-2}+\vep_{n-1}-\vep_n)\qquad (j=n-1),&\\
&&&\phantom{\bLa_j}{} =\frac{1}{2}
(\vep_1+\cdots+\vep_{n-2}+\vep_{n-1}+\vep_n)\qquad (j=n),&
\end{alignat*}
Now for $\mu\in J$, let us def\/ine $\hat{\mu}\in \cA$ by
\begin{gather*}
\hat{\mu}=\vep_{\mu}-\vep\qquad (1\leq \mu\leq n+1)\quad {\rm for }\ A_n,\\
\phantom{\hat{\mu}}{}=\pm\vep_j\ {\rm or}\ 0\qquad  (\mu=\pm j\ (1\leq j\leq n),\ {\rm or}\ \mu=0)\quad {\rm for }\ B_n,\\
\phantom{\hat{\mu}}{}=\pm\vep_j\qquad  (\mu=\pm j\ (1\leq j\leq
n))\quad {\rm for }\ C_n, D_n.
\end{gather*}
We then def\/ine a dynamical variable $a\in \h^*$ of the face
model of type $\g$ as follows.
\begin{gather*}
a_{\mu}=(a+\rho|\hat{\mu})\qquad (\mu\not=0),\\
\phantom{a_{\mu}}{}=-\frac{1}{2}\qquad (\mu=0).
\end{gather*}
We also set $a_{\mu\nu}=a_\mu-a_\nu$.
\begin{prop}\lb{dynamicalvar}
If one parameterizes $a\in \h^*$ such that
$a+\rho=\sum\limits_{i=0}^n s_i \Lambda_i$, one has
\begin{gather*}
a_{\mu}=\frac{1}{n+1}\left(-\sum_{j=1}^{\mu-1}js_j+\sum_{j=\mu}^n(n+1-j)s_j
\right) \qquad (1\leq \mu\leq n+1)\quad {\rm for }\ A^{(1)}_n,\\
\phantom{a_{\mu}}{} =\pm\left(\sum_{j=i}^{n-1}s_j
+\frac{s_n}{2}\right) \ {\rm or}\ -\frac{1}{2}
\qquad ( \mu=\pm i\ ( 1\leq i\leq n),\ {\rm or}\ 0) \quad {\rm for }\ B^{(1)}_n,\\
\phantom{a_{\mu}}{}=\pm\sum_{j=i}^n s_j \qquad ( \mu=\pm i\ ( 1\leq i\leq n) ) \quad {\rm for }\ C^{(1)}_n,\\
\phantom{a_{\mu}}{}= \pm\left(\sum_{j=i}^{n-1}
s_j+\frac{s_n-s_{n-1}}{2}\right) \quad ( \mu=\pm i\ ( 1\leq i\leq
n) ) \quad {\rm for }\ D^{(1)}_n.
\end{gather*}
\end{prop}

Note that $a_{\mu\nu}=\sum\limits_{j=\mu}^{\nu-1}s_j$ for
$A^{(1)}_n$
 and $a_{-\mu}=-a_\mu$ for $B^{(1)}_n$, $C^{(1)}_n$, $D^{(1)}_n$.

Then Jimbo--Miwa--Okado's solutions to the face type YBE are given
as follows.
\begin{gather}
W\BW{a}{b}{c}{d}{u}=\kappa(u)
\overline{W}\BW{a}{b}{c}{d}{u},\lb{JMOweight}
\\
(I)\ \ \overline{W}\BW{a}{a+\hat{\mu}}{a+\hat{\mu}}{a+2\hat{\mu}}{u}=1\quad (\mu\not=0),\nonumber\\
\phantom{(I)\ \ }
\overline{W}\BW{a}{a+\hat{\mu}}{a+\hat{\mu}}{a+\hat{\mu}+\hat{\nu}}{u}
=\frac{[1][a_{\mu\nu}-u]}{[1+u][a_{\mu\nu}]}\qquad (\mu\not=\pm\nu),\nonumber\\
\phantom{(I)\ \ }
\overline{W}\BW{a}{a+\hat{\nu}}{a+\hat{\mu}}{a+\hat{\mu}+\hat{\nu}}{u}
=\frac{[u]\sqrt{[a_{\mu\nu}+1][a_{\mu\nu}-1]}}{[1+u][a_{\mu\nu}]}\quad (\mu\not=\pm\nu),\nonumber\\
\phantom{(I)\ \ }\hspace{5cm} {\rm for}\ A_n^{(1)}, B_n^{(1)}, C_n^{(1)}, D_n^{(1)},\nonumber\\
(II)\
\overline{W}\BW{a}{a+\hat{\nu}}{a+\hat{\mu}}{a}{u}=\frac{[u][1]
[a_{\mu,-\nu}+1+\eta-u]}{[\eta-u][1+u][a_{\mu,-\nu}+1]}\sqrt{G_{a\mu}G_{a\nu}}\quad
(\mu\not=\nu),\nonumber\\
\phantom{(II)\ }{}
\overline{W}\BW{a}{a+\hat{\mu}}{a+\hat{\mu}}{a}{u}=
\frac{[\eta+u][1][a_{\mu,-\mu}+1+2\eta-u]}{[\eta-u][1+u][a_{\mu,-\mu}+1+2\eta]}\nonumber\\
\phantom{(II)\ }\qquad
{}-\frac{[u][1][a_{\mu,-\mu}+1+\eta-u]}{[\eta-u][1+u][a_{\mu,-\mu}+1+2\eta]}
\sum_{\kappa\not=\mu}\frac{[a_{\mu,-\kappa}+1+2\eta]}{[a_{\mu,-\kappa}+1]}G_{a\mu},\quad
 {\rm for}\ B_n^{(1)}, C_n^{(1)}, D_n^{(1)},\nonumber
\end{gather}
where $\eta=-{th^\vee}/{2}\ (t={({\rm long\ root})^2}/{2}$) is the
crossing parameter, and the symbol $[u]$ denotes the Jacobi
elliptic theta function
\begin{gather}
[u]= q^{r/4}e^{\pi i/4}\left(-\frac{2\pi i}{\log p}\right)^{-1/2}
q^{\frac{u^2}{r}-u}\Theta_{p}(q^{2u}),
\qquad p=q^{2r},\lb{thetafunc}\\
\Theta_{p}(z)=(z;p)_\infty(p/z;p)_\infty(p;p)_\infty,\qquad
(z;p)_\infty=\prod_{n=0}^\infty(1-zp^n).\nn
\end{gather}
The $\kappa(u)$ denotes a function satisfying the following
relations
\begin{gather*}
\kappa(u)\kappa(-u)=1\qquad {\rm for}\ A_n^{(1)}, B_n^{(1)}, C_n^{(1)}, D_n^{(1)},\\ 
\kappa(\eta-u)\kappa(\eta+u)=\frac{[1+\eta+u][1+\eta-u]}{[\eta+u][\eta-u]}\qquad {\rm for}\ A_n^{(1)},\\ 
\kappa(u)\kappa(\eta+u)=\frac{[-u][1+\eta+u]}{[1-u][\eta+u]}
\qquad {\rm for}\ B_n^{(1)}, C_n^{(1)}, D_n^{(1)}.
\end{gather*}
The $G_{a\mu}=\frac{G_{a+\hat{\mu}}}{G_a}$ denotes a ratio of the
principally specialized character $G_a$ for the dual af\/f\/ine
Lie algebra $\g^\vee$~\cite{JMO}
\begin{gather*}
G_a=\prod_{1\leq i<j\leq n+1}[a_i-a_j]\qquad {\rm for }\ A^{(1)}_n,\\
\phantom{G_a}{}=\vep(a)\prod_{i=1}^nh(a_i)\prod_{1\leq i<j\leq
n}[a_i-a_j][a_i+a_j]\qquad {\rm for }\ B^{(1)}_n, C^{(1)}_n,
D^{(1)}_n.
\end{gather*}
Here $\vep(a)$ denotes a sign factor such that
$\vep(a+\hat{\mu})/\vep(a)=s$. $s$ and $h(a)$ are listed in the
following Table
$$
\begin{array}{c|cccc}
&A_n^{(1)}&B_n^{(1)}&C_n^{(1)}&D_n^{(1)}\\ \hline
h^\vee&n+1 &2n-1 &n+1 &2n-2\\
t&1&1&2&1\\
s&1&1&-1&1\\
h(a)&1&[a]&[2a]&1\\
\end{array}
$$

\noindent {\bf Remark.} Our normalization of the weights $W$ and
some notations are dif\/ferent from those in~\cite{JMO}. Their
relations are given as follows
\begin{gather*}
W\BW{a}{b}{c}{d}{u}=\frac{[1]}{[1+u]}\kappa(u)W_{JMO}\BW{a}{b}{c}{d}{u}
\qquad {\rm for }\
A^{(1)}_n,\\
\phantom{W\BW{a}{b}{c}{d}{u}}{}=\frac{[1][\eta]}{[1+u][\eta-u]}\kappa(u)W_{JMO}\BW{a}{b}{c}{d}{u}
\qquad {\rm for }\
B^{(1)}_n, C^{(1)}_n, D^{(1)}_n,\\
r = L_{JMO},\qquad \log p = \frac{4\pi^2}{\log p_{JMO}}.
\end{gather*}
Here the symbols with subindex ${JMO}$ denote the ones
in~\cite{JMO}. Note $[ u ] =\left[{u}\right]_{JMO}$.

The following theorem states basic properties of the face weights.
\begin{thm}[\cite{JMO}]
The face weight $W$ satisfies i) the face type Yang--Baxter
equation, ii)~the first and iii) the second inversion relations
\begin{alignat*}{3}
& i)\ \ \ &&\sum_g W\BW{f}{g}{e}{d}{u}W\BW{a}{b}{f}{g}{u+v}W\BW{b}{c}{g}{d}{v}&\\
&&&\quad{}=\sum_gW\BW{a}{b}{g}{c}{u}W\BW{g}{c}{e}{d}{u+v}W\BW{a}{g}{f}{e}{v},&\\
& ii)&&\sum_{g}W\BW{a}{g}{d}{c}{u}W\BW{a}{b}{g}{c}{-u}=\delta_{bd},&\\
&
iii)&&\sum_{g}\frac{G_aG_g}{G_bG_d}W\BW{a}{b}{d}{g}{-u}W\BW{c}{d}{b}{g}{2\eta+u}
=\delta_{ac}.&
\end{alignat*}
In addition, we have the crossing symmetry except for
$\g=A_n^{(1)}$ $(n>1)$
\begin{gather*}
iv)\ \ W\BW{a}{b}{c}{d}{u}=\sqrt{\frac{G_b G_c}{G_aG_d}}
W\BW{c}{a}{d}{b}{\eta-u}.
\end{gather*}
\end{thm}

The following theorem is communicated by Jimbo and Okado and is
not written explicitly in~\cite{JMO}.
\begin{thm}\lb{partII}
For $\g=B^{(1)}_n, C^{(1)}_n, D^{(1)}_n$, the weights listed in
the part (II) of \eqref{JMOweight} is determined uniquely from
those in  (I) by requiring the face type Yang--Baxter equation and
the crossing symmetry relations.
\end{thm}

\begin{figure}[h]
\centering
\includegraphics[width=16cm]{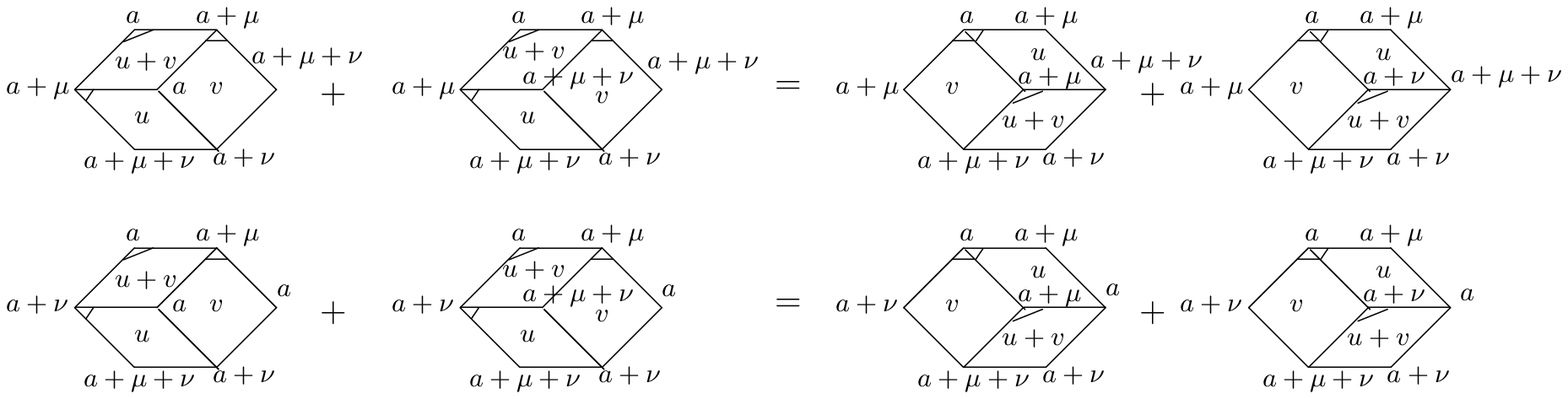}
\caption{Two relevant equations ($\mu\not=\pm \nu$).}
\label{Konno-figure1}
\end{figure}

\begin{proof}[Sketch of proof.] It is easy to see that the f\/irst type of weights in (II)
is determined by those in (I) by using the crossing symmetry
relation. Then the second type of weights in (II) is determined by
solving the system of two linear equations (YBE)  shown in
Fig.~\ref{Konno-figure1}. Here unknowns are the weights
$W\mat{{a}&{a+\mu}\cr {a+\mu}&{a}\cr}$ and
$W\mat{{a+\nu}&{a+\mu+\nu}\cr {a+\mu+\nu}&{a+\nu}\cr}$, and the
other weights are in (I).
\end{proof}

\subsection[The difference equations for the twistor]{The dif\/ference equations for the twistor}
We here solve the dif\/ference equation for the twistor. Then
using Theorem \ref{RmatCmat}, we derive the dynamical $R$ matrix
$R_{VV}(z,\la)$ as the connection matrix in the vector
representation $(\pi_V, V)$, and argue that it coincides with
Jimbo--Miwa--Okado's solution up to a gauge transformation.

Let us consider the dif\/ference equation \eqref{Fdiff} in the
vector representation.
\begin{gather}
F(pz,\la)=(q^{2\pi_V(\bar{\theta}(\lambda))} \otimes
\id)F(z,\la)(q^{-2\pi_V(\bar{\theta}(\lambda))}\otimes
\id)q^{\pi_{V\otimes V}(T)}R(pz),\lb{Fdiffvec}
\end{gather}
where $\la$ is parameterized as \eqref{parametlamda},
$\bar{\theta}(\la)$ is given by \eqref{bartheta}, $T=c\otimes
\Lambda_0+\Lambda_0\otimes c +\sum\limits_{i=1}^n\bar{h}_i\otimes
\bar{h}^i$, $F(z,\la)=(\pi_V\otimes \pi_V)F(z,\la)$, and
$R(z)=(\pi_V\otimes \pi_V)\cR(z)$.

Let $\{v_j\ | j\in J\}$ be a basis of $V$ and $E_{i,j}$ be the
matrix unit def\/ined by $E_{i,j}v_k=\delta_{j,k}v_i$. The action
of the generators on $V$ is given by \cite{DO,DavO,JinMisO} (for
$C_n^{(1)}$, the conventions used here are slightly dif\/ferent
from \cite{JinMisO})
\begin{gather*}
\pi_V(e_0)=E_{n+1,1} \qquad {\rm for}\ A^{(1)}_n,\\
\phantom{\pi_V(e_0)}{}=(-)^n(E_{-1,2}-E_{-2,1}) \qquad {\rm for}\ B^{(1)}_n,\\
\phantom{\pi_V(e_0)}{}=(-)^{n-1}(E_{-1,2}-E_{-2,1}) \qquad {\rm for}\ D^{(1)}_n,\\
\phantom{\pi_V(e_0)}{}=E_{-1,1} \qquad {\rm for}\ C^{(1)}_n,\\
\pi_V(e_i)=E_{i,i+1} \qquad (1\leq i\leq n)\quad {\rm for}\ A^{(1)}_n,\\
\phantom{\pi_V(e_i)}{}=E_{i,i+1}-E_{-i-1,-i}\qquad (1\leq i\leq n-1)\quad {\rm for}\ B^{(1)}_n, D^{(1)}_n, C^{(1)}_n,\\
\pi_V(e_n)=\sqrt{[2]_{q_n}}(E_{n,0}-E_{0,-n})\qquad {\rm for}\ B^{(1)}_n,\\
\phantom{\pi_V(e_n)}{}=E_{n,-n}\qquad {\rm for}\ C^{(1)}_n,\\
\phantom{\pi_V(e_n)}{} =E_{n-1,-n}-E_{n,-n+1} \qquad {\rm for}\
D^{(1)}_n,
\\
\pi_V(t_0)=\sum_{j\in J} q^{-\delta_{j,1}+\delta_{j,n+1}}E_{j,j} \qquad {\rm for}\ A^{(1)}_n,\\
\phantom{\pi_V(t_0)}{}=\sum_{j\in J} q^{-\delta_{j,1}-\delta_{j,2}+\delta_{j,-1}+\delta_{j,-2}}E_{j,j} \qquad {\rm for}\ B^{(1)}_n, D_n^{(1)}\\
\phantom{\pi_V(t_0)}{}=\sum_{j\in J} q^{-2\delta_{j,1}+2\delta_{j,-1}}E_{j,j} \qquad {\rm for}\ C^{(1)}_n,\\
\pi_V(t_i)=\sum_{j\in J} q^{\delta_{j,i}-\delta_{j,i+1}}E_{j,j}\qquad (1\leq i\leq n)\quad {\rm for}\ A^{(1)}_n\\
\phantom{\pi_V(t_i)}{}=\sum_{j\in J}
q^{\delta_{j,i}-\delta_{j,i+1}
+\delta_{j,-i-1}-\delta_{j,-i}}E_{j,j}\qquad (1\leq i\leq
n-1)\quad {\rm for}\ B^{(1)}_n, C^{(1)}_n, D^{(1)}_n,
\\
\pi_V(t_n)=\sum_{j\in J} q^{\delta_{j,n}-\delta_{j,-n}}E_{j,j} \qquad {\rm for}\ B^{(1)}_n,\\
\phantom{\pi_V(t_n)}{}
=\sum_{j\in J} q^{2\delta_{j,n}-2\delta_{j,-n}}E_{j,j}\qquad {\rm for}\ C^{(1)}_n,\\
\phantom{\pi_V(t_n)}{} =\sum_{j\in J}
q^{\delta_{j,n-1}+\delta_{j,n}-\delta_{j,-n}-\delta_{j,-n+1}}E_{j,j}
\qquad {\rm for}\ D^{(1)}_n,
\end{gather*}
and  $\pi_V(f_i)=\pi_V(e_i)^t$.

A basis $\{\bar{h}_i\}$ of $\bar{\h}$ and its dual basis
$\{\bar{h}^i\}$ w.r.t $(\cdot |\cdot )$ are given as follows
\begin{alignat*}{3}
& A_n:\quad &&  \pi_V(\bar{h}_i)=E_{i,i}-E_{i+1,i+1}\qquad (1\leq i\leq n),& \\
&&& \pi_V(\bar{h}^i)=\frac{1}{n+1}\left(
(n-i+1)\sum_{j=1}^iE_{j,j}-i\sum_{j=i+1}^{n+1}E_{j,j}\right)
\qquad (1\leq i\leq n),& \\
& B_n :\quad && \pi_V(\bar{h}_i)=E_{i,i}-E_{i+1,i+1}+E_{-i-1,-i-1}-E_{-i,-i}\qquad (1\leq i\leq n-1), &\\
&&& \pi_V(\bar{h}_n)=2(E_{n,n}-E_{-n,-n}),&\\
&&& \pi_V(\bar{h}^i)=\sum_{j=1}^i(E_{j,j}-E_{-j,-j}) \qquad (1\leq i\leq n-1),&\\
&&& \pi_V(\bar{h}^n)=\frac{1}{2}\sum_{j=1}^n(E_{j,j}-E_{-j,-j}),&\\
& C_n :\quad && \pi_V(\bar{h}_i)=E_{i,i}-E_{i+1,i+1}+E_{-i-1,-i-1}-E_{-i,-i}\qquad (1\leq i\leq n-1),&\\
&&& \pi_V(\bar{h}_n)=E_{n,n}-E_{-n,-n},&\\
&&& \pi_V(\bar{h}^i)=\sum_{j=1}^i(E_{j,j}-E_{-j,-j}) \qquad (1\leq i\leq n),&\\
& D_n :\quad && \pi_V(\bar{h}_i)=E_{i,i}-E_{i+1,i+1}+E_{-i-1,-i-1}-E_{-i,-i}\qquad (1\leq i\leq n-1),&\\
&&& \pi_V(\bar{h}_n)=E_{n-1,n-1}+E_{n,n}-E_{-n,-n}-E_{-n+1,-n+1},&\\
&&& \pi_V(\bar{h}^i)=\sum_{j=1}^i(E_{j,j}-E_{-j,-j}) \qquad (1\leq i\leq n-2),&\\
&&& \pi_V(\bar{h}^{n-1})=\frac{1}{2}\sum_{j=1}^{n-1}(E_{j,j}-E_{-j,-j})- \frac{1}{2}(E_{n,n}-E_{-n,-n}),&\\
&&&
\pi_V(\bar{h}^{n})=\frac{1}{2}\sum_{j=1}^{n-1}(E_{j,j}-E_{-j,-j})+
\frac{1}{2}(E_{n,n}-E_{-n,-n}).&
\end{alignat*}
Then one can easily verify the following.
\begin{prop}\lb{qT}
\begin{gather*}
q^{\pi_{V\otimes V}(T)}=
q^{-\frac{1}{n+1}}\sum_{i,j\in J}q^{\delta_{i,j}}E_{i,i}\otimes E_{j,j}\qquad {\rm for }\ A^{(1)}_n,\\
\phantom{q^{\pi_{V\otimes V}(T)}}{} =\sum_{i,j\in
J}q^{\delta_{i,j}-\delta_{i,-j}}E_{i,i}\otimes E_{j,j}\qquad {\rm
for }\ B^{(1)}_n, C^{(1)}_n, D^{(1)}_n.
\end{gather*}
\end{prop}

\begin{prop}\lb{qbtheta}
If we parameterize $\bar{\la}$ such that
 $\bar{\la}=\sum\limits_{i=1}^n (s_i+1)\bar{h}^i$, we have
\begin{gather*}
q^{-2\pi_V(\bar{\theta}(\la))} = q^{\frac{n}{n+1}}\sum_{j\in J}
q^{2a_j}E_{j,j}
\qquad {\rm for }\ A^{(1)}_n,\\
\phantom{q^{-2\pi_V(\bar{\theta}(\la))}}{} =\sum_{j\in
J}q^{2a_j+1}E_{j,j}\qquad {\rm for }\ B^{(1)}_n, C^{(1)}_n,
D^{(1)}_n,
\end{gather*}
where $a_j$ $(j\in J)$ is given by Proposition {\rm
\ref{dynamicalvar}}.
\end{prop}

The $R$ matrix $R(z)$ of $U_q(\g)$ in the vector representation is
well known \cite{Jimbo86,DO,DavO,JinMisO} (for~$C^{(1)}_n$, we
modif\/ied the $R$ matrix in \cite{JinMisO} according to the
convention used here)
\begin{gather}
R(z)=\rho(z)\left\{
\sum_{i\in J\atop i\not=0}E_{i,i}\otimes E_{i,i}+b(z)\sum_{i,j\atop i\not=\pm j}E_{i,i}\otimes E_{j,j}\right. \nn\\
\phantom{R(z)=}{} +\sum_{i\prec j\atop i\not= -j}\Bigl(c(z)E_{i,j}\otimes E_{j,i}+zc(z)E_{j,i}\otimes E_{i,j}\Bigr)\nn\\
\left.\phantom{R(z)=}{} +\frac{1}{(1-q^2z)(1-\xi
z)}\sum_{i,j}a_{ij}(z)E_{i,j}\otimes E_{-i,-j}
\right\},\lb{trigRmat}\\
b(z)=\frac{q(1-z)}{1-q^2z}, \quad c(z)=\frac{1-q^2}{1-q^2z},\nn\\
\rho(z)=q^{-\frac{n}{n+1}}\frac{(q^2z;\xi^2)_\infty(q^{-2}\xi^2z;\xi^2)_\infty}
{(z;\xi^2)_\infty(\xi^{2}z;\xi^2)_\infty}\qquad {\rm for}\ A^{(1)}_n,\nn\\
\phantom{\rho(z)}{}=q^{-1}\frac{(q^2z;\xi^2)_{\infty}(\xi
z;\xi^2)^2_{\infty}(q^{-2}\xi^2 z;\xi^2)_{\infty}}
{(z;\xi^2)_{\infty}(q^{-2}\xi z;\xi^2)_{\infty}(q^2\xi
z;\xi^2)_{\infty}(\xi^2z;\xi^2)_{\infty}} \qquad {\rm for}\
B^{(1)}_n, C^{(1)}_n, D^{(1)}_n,\nn
\\
a_{ij}(z)=0 \qquad {\rm for}\ A^{(1)}_n,\nn\\
\phantom{a_{ij}(z)}{}=\left\{\begin{array}{ll}
(q^2-\xi z)(1-z)+\delta_{i,0}(1-q)(q+z)(1-\xi z)\quad & (i=j),\\
(1-q^2)[\e_i\e_jq^{\bar{j}-\bar{i}}(z-1)+\delta_{i,-j}(1-\xi z)]& (i\prec j),\\
(1-q^2)z[\xi \e_i\e_j q^{\bar{j}-\bar{i}}(z-1)+\delta_{i,-j}(1-\xi
z)]& (i\succ j),
\end{array}\right. \quad {\rm for}\ B^{(1)}_n, C^{(1)}_n, D^{(1)}_n. \nn
\end{gather}
Here $\xi=q^{t h^\vee}$, and $\e_j=1$ $(j>0$), $-1$ $(j<0)$ for
$\g=C_n^{(1)}$ and $\e_j=1$ $(j\in J)$ for the other cases. The
symbol $\bar{j}$ is def\/ined by
\begin{gather*}
\bar{j}=\left\{\begin{array}{ll} j-\e_j&(j=1,\ldots, n),\\
                          n-\e_j& (j=0), \\
                          j+N-\e_j\quad& (j=-n,\ldots,-1).\end{array}\right.
\end{gather*}
Then due to the formula \eqref{facetwistor}, we make the following
ansatz for the twistor $F(z,\la)$ in the vector representation.
\begin{gather}
F(z,\la)=f(z)\left\{ \sum_{i\in J\atop i\not=0}E_{i,i}\otimes
E_{i,i}+\sum_{i,j\atop i\not=\pm j}X_{ij}^{ij}(z)
E_{i,i}\otimes E_{j,j}\right. \lb{Fmat}\\
\phantom{F(z,\la)=}{} +\sum_{i\prec j\atop i\not=
-j}\Bigl(X_{ij}^{ji}(z)E_{i,j}\otimes E_{j,i}+
X_{ji}^{ij}(z)E_{j,i}\otimes E_{i,j}\Bigr) \left.
+\sum_{i,j}X_{i,-i}^{j,-j}(z)E_{i,j}\otimes E_{-i,-j} \right\},\nn
\end{gather}
where $X_{ij}^{kl}$ denote  unknown functions to be determined.

From the from of $R(z)$ and $F(z,\la)$ in \eqref{trigRmat} and
\eqref{Fmat}, one f\/inds that the dif\/ference
equation~\eqref{Fdiffvec} consists of $1\times 1$, $2\times 2$ and
$N\times N$ blocks. The numbers of blocks of each size contained
in the equation are listed as follows
$$
\begin{array}{c|ccc}
          &1\times 1&2\times2&N\times N\\ \hline
A^{(1)}_n &n+1&\frac{n(n+1)}{2}&0\\
B^{(1)}_n &2n &2n^2&1\\
C^{(1)}_n &2n&2n(n-1)&1\\
D^{(1)}_n &2n &2n(n-1)&1\\
\end{array}
$$
By using Propositions \ref{qT}, \ref{qbtheta} and
\eqref{trigRmat}, \eqref{Fmat},
 we obtain the following equations.

 $1\times 1$ blocks:
\begin{gather*}
f(pz)=q^{\frac{n}{n+1}}\rho(pz)f(z)\qquad {\rm for}\ A_n,\\
\phantom{f(pz)}{}=q \rho(pz)f(z)\qquad {\rm for}\ B_n, C_n, D_n.
\end{gather*}

$2\times 2$ blocks:
\begin{gather*}
\mat{X_{ij}^{ij}(pz)&X_{ij}^{ji}(pz)\cr
     X_{ji}^{ij}(pz)&X_{ji}^{ji}(pz)\cr}=q^{-1}\mat{X_{ij}^{ij}(z)&w_{ij}^{-1}X_{ij}^{ji}(z)\cr
     w_{ij}X_{ji}^{ij}(z)&X_{ji}^{ji}(z)\cr}
\mat{b(pz)&c(pz)\cr
     pz c(pz)&b(pz)\cr},\\
\phantom{\mat{X_{ij}^{ij}(pz)&X_{ij}^{ji}(pz)\cr
     X_{ji}^{ij}(pz)&X_{ji}^{ji}(pz)\cr}=}{}
 (i,j\in J,\ i\prec j,\ i\not=-j)
\end{gather*}
where we set $w_{ij}=q^{2(a_{i}-a_j)}$.

$N\times N$ block:
\begin{gather*}
X_{i,-i}^{j,-j}(pz) =\frac{q^{-2}}{(1-pq^2z)(1-p\xi z)}\sum_{k\in
J} q^{-2(a_i-a_k)}a_{kj}(pz) X_{i,-i}^{k,-k}(z)\qquad (i,j\in J).
\end{gather*}
Here  we dropped a scalar factor in the $2\times 2$ and $N\times
N$ blocks by using the equation in the $1\times 1$ block.

Note that the dif\/ference equations in the $2\times 2$ blocks
have the same structure as the one in  the case $\g=\slth$, which
was analyzed completely in~\cite{JKOS}. Let us summarize the
essence of it. The $2\times 2$ block equation consists of two
 2nd order $q$-dif\/ference equations of the type
\[
(q^c-q^{a+b+1}z)u(q^2z)-\{(q+q^{c})-(q^a+q^b)qz\}u(qz)+q(1-z)u(z)=0.
\]
This  equation has two independent solutions of the form
$z^\alpha\sum\limits_{n=0}^\infty a_nz^n$ around $z=0$, which are
given by the basic hypergeometric series
\[
{}_2\phi_1\left(\mmatrix{q^a\quad q^b\cr
                               q^c\cr };q,z\right)
                               =\sum_{n=0}^\infty\frac{(q^a;q)_n(q^b;q)_n}{(q^c;q)_n(q;q)_n}z^n,
\]
and \be z^{1-c}{}_2\phi_1\left(\mmatrix{q^{a-c+1}\quad
q^{b-c+1}\cr
                               q^{2-c}\cr};q,z\right),
\en where $(x;q)_n=\prod\limits_{j=0}^{n-1}(1-x q^j)$,
$(x;q)_0=1$. The connection formula for these solutions is well
known:
\begin{gather}
{}_2\phi_1\left(\mmatrix{q^a\quad q^b\cr
                               q^c\cr };q,1/z\right)=
                               \frac{\Gamma_q(c)\Gamma_q(b-a)\Theta_q(q^{1-a}z)}{\Gamma_q(b)\Gamma_q(c-a)\Theta_q(qz)}
{}_2\phi_1\left(\mmatrix{q^a\quad q^{a-c+1}\cr
                               q^{a-b+1}\cr };q,q^{c-a-b+1}z\right)\nn\\
\phantom{{}_2\phi_1\left(\mmatrix{q^a\quad q^b\cr
                               q^c\cr };q,1/z\right)=}{}
+\frac{\Gamma_q(c)\Gamma_q(a-b)\Theta_q(q^{1-b}z)}{\Gamma_q(a)\Gamma_q(c-b)\Theta_q(qz)}
{}_2\phi_1\left(\mmatrix{q^b\quad q^{b-c+1}\cr
                               q^{b-a+1}\cr };q,q^{c-a-b+1}z\right),\lb{conn}
\end{gather}
where
\[
\Gamma_q(z)=\frac{(q;q)_\infty}{(q^z;q)_\infty}(1-q)^{1-z}.
\]
By using this, one can derive the connection matrices for the
$2\times 2$ block parts.

In our case, the initial condition \eqref{Fini} leads to
\begin{gather*}
f(0)=1,\\
\mat{X_{ij}^{ij}(0)&X_{ij}^{ji}(0)\cr
     X_{ji}^{ij}(0)&X_{ji}^{ji}(0)\cr}=\mat{1&\frac{(q-q^{-1})w_{ij}}{1-w_{ij}}\cr
                                              0&1  }.
\end{gather*}
The solutions to the $1\times 1$ and $2\times 2$ blocks are given
as follows.

$1\times 1$ block:
\begin{gather*}
f(z)=\frac{\{pz\}\{p\xi^2z\}}
{\{pq^2z\}\{pq^{-2}\xi^{2}z\}}\qquad {\rm for}\ A_n,\\
\phantom{f(z)}{}=\frac{\{pz\}\{pq^{-2}\xi z\}\{pq^2\xi
z\}\{p\xi^2z\}}{\{pq^2z\}\{p\xi z\}^2 \{pq^{-2}\xi^2 z\}} \qquad
{\rm for}\ B_n, C_n, D_n,
\end{gather*}
where
\[
\{z\}=\prod_{n,m=0}^\infty(1-z\xi^{2n}p^m).
\]

$2\times 2$ block:
\begin{gather*}
X_{ij}^{ij}(z)={}_2\phi_1\left(\mmatrix{w_{ij}q^2\quad q^2\cr
                                                   w_{ij}\cr };p,pq^{-2}z\right),\\
X_{ij}^{ji}(z)=\frac{(q-q^{-1})w_{ij}}{1-w_{ij}}
{}_2\phi_1\left(\mmatrix{w_{ij}q^2\quad pq^2\cr
                                                   pw_{ij}\cr };p,pq^{-2}z\right),\\
X_{ji}^{ij}(z)=\frac{(q-q^{-1})pw^{-1}_{ij}}{1-pw^{-1}_{ij}z}
{}_2\phi_1\left(\mmatrix{pw^{-1}_{ij}q^2\quad pq^2\cr
                                                   p^2w^{-1}_{ij} };p,pq^{-2}z\right),\\
X_{ji}^{ji}(z)={}_2\phi_1\left(\mmatrix{pw^{-1}_{ij}q^2\quad
q^2\cr
                                                   pw^{-1}_{ij} };p,pq^{-2}z\right).
\end{gather*}
Then due to the formulae \eqref{conn} and \eqref{DynamicalR} or
Theorem \ref{RmatCmat}, we determine the $1\times 1$ and $2\times
2$ blocks of the dynamical $R$ matrix
\begin{gather*}
(\pi_V\otimes \pi_V)\cR(z,\la)\\
\quad{}=\rho_{ell}(z)\left\{
\sum_{i\in J\atop i\not=0}E_{i,i}\otimes E_{i,i}+\sum_{i\prec j\atop i\not=\pm j}\Bigl(R_{ij}^{ij}(z,w_{ij})E_{i,i}\otimes E_{j,j}+R_{ji}^{ji}(z,w_{ij})E_{j,j}\otimes E_{i,i}\Bigr)\right. \\
\qquad {}+\sum_{i\prec j\atop i\not= -j}\Bigl(R_{ij}^{ji}(z,w_{ij})E_{i,j}\otimes E_{j,i}+R^{ij}_{ji}(z,w_{ij})E_{j,i}\otimes E_{i,j}\Bigr)\\
\qquad \left. {} +\sum_{i,j}R_{i-i}^{j-j}(z,w_{ij})E_{i,j} \otimes
E_{-i,-j} \right\}
\end{gather*}
as follows.

$1\times 1$ block:
\begin{gather*}
\rho_{ell}(z)={f(z^{-1})\rho(z)}{f(z)^{-1}}\\
\phantom{\rho_{ell}(z)}{}=q^{-\frac{n}{n+1}}\frac{\{q^2z\}\{q^{-2}\xi^2z\}\{p/z\}\{p\xi^2/z\}}
{\{z\}\{\xi^2z\}\{pq^2/z\}\{pq^{-2}\xi^2/z\}}\qquad {\rm for}\ A_n^{(1)},\\
\phantom{\rho_{ell}(z)}{}=q^{-1}\frac{\{q^2z\}\{\xi
z\}^2\{q^{-2}\xi^2z\}\{p/z\}\{pq^{-2}\xi/z\}
\{pq^{2}\xi/z\}\{p\xi^2/z\}}
{\{z\}\{q^{-2}\xi z\}\{q^{2}\xi z\}\{\xi^2z\}\{pq^2/z\}\{p\xi/z\}^2\{pq^{-2}\xi^2/z\}}\\
\phantom{\rho_{ell}(z)=}{} {\rm for}\ B_n^{(1)},\ C_n^{(1)},\
D_n^{(1)}.
\end{gather*}

$2\times 2$ blocks: for $i\prec j$, $i\not= -j$,
\begin{gather*}
R_{ij}^{ij}(z,w_{ij})=q\frac{(pw_{ij}^{-1}q^{2};p)_\infty(pw_{ij}^{-1}q^{-2};p)_\infty}
{(pw_{ij}^{-1};p)^2_\infty}\frac{\Theta_p(z)}{\Theta_p(q^2z)},\\
R_{ji}^{ji}(z,w_{ij})=q\frac{(w_{ij}q^{2};p)_\infty(w_{ij}q^{-2};p)_\infty}
{(w_{ij};p)^2_\infty}\frac{\Theta_p(z)}{\Theta_p(q^2z)},\\
R_{ij}^{ji}(z,w_{ij})=\frac{\Theta_p(q^2)}{\Theta_p(w_{ij})}\frac{\Theta_p(w_{ij}z)}{\Theta_p(q^2z)},\\
R^{ij}_{ji}(z,w_{ij})=z\frac{\Theta_p(q^2)}{\Theta_p(pw^{-1}_{ij})}\frac{\Theta_p(pw^{-1}_{ij}z)}{\Theta_p(q^2z)}.
\end{gather*}
By setting $z=q^{2u}$ and using \eqref{thetafunc}, we can
reexpress these matrix elements in terms of the theta functions
with some extra factors including $q$ with fractional power and
inf\/inite products. Then making an appropriate
 gauge transformation, we can sweep away all the extra factors and
f\/ind that the $1\times 1$ and $2\times 2$ block parts coincide
with the part (I) of Jimbo--Miwa--Okado's solution, i.e.
\[
R_{ij}^{kl}(z,w_{ij})\Leftrightarrow
{W}\BW{a}{a+\hat{k}}{a+\hat{i}}{a+\hat{i}+\hat{j}}{z} \qquad {\rm
in}\ (I) .
\]
From 2) of Theorem \ref{connectMat}, \eqref{Gcross} and Theorem
\ref{partII}, the remaining part (II) is determined uniquely from
the part (I).  We hence obtain the following theorem.
\begin{thm}\lb{VecRep}
For $\g=A_n^{(1)}, B_n^{(1)}, C_n^{(1)}, D_n^{(1)}$, the vector
representation of the universal dyna\-mical $R$ matrix $\cR(\la)$
coincides with Jimbo--Miwa--Okado's elliptic solutions to the face
type YBE.
\end{thm}
To solve the dif\/ference equation in the $N\times N$ block
directly is an open problem.

For the cases $\g$ being the twisted af\/f\/ine Lie algebras
$A^{(2)}_{2n}$ and $A^{(2)}_{2n-1}$, Kuniba derived elliptic
solutions to the face type YBE. His construction is based on a
common structure of the $R$ matrices of the twisted $U_q(\g)$
 to those of the  $B_n^{(1)}$, $C_n^{(1)}$, $D_n^{(1)}$ types.
 In fact, the resultant face weights have the common
 $2\times 2$ block part, as a function of $a_\mu$, to the cases $\g=A^{(1)}_n, B_n^{(1)}, C_n^{(1)}, D_n^{(1)}$.
The simplest $A_2^{(2)}$ case was investigated in \cite{KK04}.
 In view of these facts, we expect the same statement
 as Theorem \ref{partII} is valid in the twisted cases, too.
\begin{conj}
Similar statement to Theorem {\rm \ref{VecRep}} is true for
Kuniba's solution of $A_{2n}^{(2)}$, $A_{2n-1}^{(2)}$ types and
for Kuniba--Suzuki's solution of $G_2^{(1)}$ type.
\end{conj}

\appendix
\section{Proof of Lemma \ref{DeqK}}
We here give a direct proof of Lemma \ref{DeqK} and leave a
derivation of the $q$-KZ Equation\eqref{qKZeq} from it as an
exercise.

Let $\cR$ be the universal $R$ matrix of $U_q(\g)$ and  write
\begin{gather}
\cR=\sum_{j}a_j\otimes b_j,\lb{Rab}
\end{gather}
and set $\cU=\sum_j S(b_j)a_j=\sum_j b_jS^{-1}(a_j)$ and
$\cZ=q^{2\rho}\cU$.

\begin{lem}[\cite{Drinfeld}]\lb{DriCas}
\begin{gather*}
(1)\ \ \cU x\cU^{-1}=S^2(x)\qquad \forall\, x\in U_q,\\
(2)\ \ \cZ x \cZ^{-1}=x, \\
(3)\ \ \cZ|_{V(\la)}=q^{(\la|\la+2\rho)}\id_{V(\la)}.
\end{gather*}
\end{lem}

\begin{lem}\lb{abb}
For \eqref{Rab},
\[
\sum_ja_j\otimes \Delta(b_{j})=\sum_{i,j}a_ia_j\otimes b_j\otimes
b_i.
\]
\end{lem}

\begin{proof}
The statement follows $(\id \otimes
\Delta)\cR=\cR^{(13)}\cR^{(12)}$.
\end{proof}

\begin{lem}\lb{aPsi}
Let $\Psi(z)$ denote a vertex operator. Then we have \be
(\id\otimes a)\Psi(z)=\sum (S(a_{(1)})\otimes 1)\Psi(z) a_{(2)},
\en where  we write $\Delta(a)=\sum a_{(1)}\otimes a_{(2)}$.
\end{lem}

\begin{proof}
\begin{gather*}
RHS =\sum (S(a_{(1)})\otimes 1)\Delta(a_{(2)})\Psi(z) \\
\phantom{RHS}=\sum  (S(a_{(1)})a_{(2)}'\otimes a_{(2)}'')\Psi(z) \\
\phantom{RHS}=\sum  (S(a_{(1)}')a_{(1)}''\otimes a_{(2)})\Psi(z) \\
\phantom{RHS}=\sum  (1\otimes \e(a_{(1)})a_{(2)})\Psi(z) \\
\phantom{RHS}=LHS.
\end{gather*}
Here we wrote $\Delta(a_{(2)})=\sum a_{(2)}'\otimes a_{(2)}''$
etc. and used $(\Delta\otimes \id)\Delta(a)=(\id\otimes
\Delta)\Delta(a)$ in the 3rd line and $m(S\otimes
\id)\Delta(a)=\e(a)$ in the 4th line.
\end{proof}

\begin{proof}[Proof of Lemma \ref{DeqK}]
Let $\la, \mu, \nu\in \h^*$ be level-$k$ elements. Let us set
$\tp=q^{2(k+h^\vee)}$ and consider
\[
\tPsi(z_1,  z_2)=\big\langle\id\otimes \id\otimes
u_\nu^*,(\id\otimes \tPsi^{\nu}_\mu(z_2)\cU) \tPsi^\mu_\la(\tp
z_1)u_\la\big\rangle,
\]
where we abbreviate $\tPsi^{\nu, v_i}_\mu(z_2)$ and
$\tPsi^{\mu,w_j}_\la(z_1)$ as $\tPsi^{\nu}_\mu(z_2)$ and
$\tPsi^\mu_\la( z_1)$, respectively. We regard $\cU$ and its
expression in terms of $a_j$, $b_j$ as  certain images of
appropriate representations of $U_q(\g)$ in the following
processes.  We evaluate $\tPsi(z_1,  z_2)$ in the following two
ways.

1) Substituting $\cU=q^{-2\rho}\cZ$ and using the intertwining
property \eqref{Psiint} and Lemma \ref{DriCas} (3), we have
\begin{gather*}
\tPsi(z_1,z_2)=(\id\otimes q^{-2\rho})q^{-(\nu|2\rho)+(\mu|\mu+2\rho)}J_{WV}(\tp z_1,z_2;\la)(w_j\otimes v_i)\\
\phantom{\tPsi(z_1,z_2)}{} =(\id\otimes q^{-2\rho})q^{-(\nu|2\rho)+(\mu|\mu+2\rho)+2(\wt(w_j)|\rho+\la)-(\wt(w_j)|\wt(w_j))}\\
\phantom{\tPsi(z_1,z_2)=}{} \times J_{WV}(\tp z_1,z_2;\la)
(q^{-2\pi_W(\bar{\theta}(\la))}\otimes \id)(w_j\otimes v_i).
\end{gather*}
In the last line, we used $(q^{-2\pi_W(\bar{\theta}(\la))}\otimes
\id)|_{w_j\otimes v_i}
=q^{-2(\wt(w_j)|\rho+\la)+(\wt(w_j)|\wt(w_j))}$.

2) Using $\cU=\sum_j b_jS^{-1}(a_j)$ and \eqref{abb}, we have
\begin{gather*}
(\id\otimes \tPsi^\nu_\mu(z_2)\cU)\tPsi^\mu_\la( z_1)
=\sum_j(\id\otimes \Delta(b_j)\tPsi^\nu_\mu(z_2)S^{-1}(a_j))\tPsi^\mu_\la( z_1)\\
\phantom(\id\otimes \tPsi^\nu_\mu(z_2)\cU)\tPsi^\mu_\la( z_1){}{}
=\sum_{i,j}(\id\otimes (b_i\otimes b_j)\tPsi^\nu_\mu(z_2)S^{-1}(a_ja_i))\tPsi^\mu_\la( z_1)\\
\phantom(\id\otimes \tPsi^\nu_\mu(z_2)\cU)\tPsi^\mu_\la( z_1){}{}
=\sum_{i,j}(\id\otimes (S(b_i)\otimes S(b_j))\tPsi^\nu_\mu(z_2)a_ia_j)\tPsi^\mu_\la( z_1)\\
\phantom(\id\otimes \tPsi^\nu_\mu(z_2)\cU)\tPsi^\mu_\la( z_1){}{}
=\sum_{i,j}(\id\otimes (S(b_i)\otimes
S(b_j))\tPsi^\nu_\mu(z_2))(1\otimes a_i)(1\otimes
a_j)\tPsi^\mu_\la( z_1)
\end{gather*}
In the 3rd line we used $(\id\otimes S)\cR=(S^{-1}\otimes
\id)\cR$. Then apply Lemma \ref{aPsi} and Lemma \ref{abb} twice
each, we have
\[
(\id\otimes \tPsi^\nu_\mu(z_2)\cU)\tPsi^\mu_\la( z_1)
=\sum_{i,j,k,l}(S(a_i)S(a_j)\otimes (S(b_ib_k)\otimes
S(b_jb_l))\tPsi^\nu_\mu(z_2))\tPsi^\mu_\la( z_1)a_ka_l.
\]
Take the expectation value $\langle\id\otimes \id\otimes u_\nu^*,
u_\la\rangle$, and use
\begin{gather*}
\left<u_\nu^*,\sum_lS(b_l) a_l u_\la\right> = q^{(\la|\nu)},\\
\sum_k S(b_k) \otimes a_k u_\la =q^{k\Lambda_0+\bar{\la}}\otimes u_\la,\\
\sum_j S(a_j) \otimes  u^*_\nu S(b_j)
=q^{-k\Lambda_0-\bar{\nu}}\otimes u_\nu^*.
\end{gather*}
Noting further that \eqref{Psiweight} implies $\tPsi_\la^\mu(\tp
z)=(\tp^{\Lambda_0}\otimes \id)\tPsi_\la^\mu(z)$, we obtain
\begin{gather*}
\tPsi(z_1,z_2)=(\tp^{\Lambda_0}\otimes \id)q^{(\la|\nu)}(1\otimes
q^{k\Lambda_0+\bar{\la}})(q^{-k\Lambda_0-\bar{\nu}}\otimes 1)\\
\phantom{\tPsi(z_1,z_2)=}{}\times
\sum_i(S(a_i)\otimes S(b_i))J_{WV}(z_1,z_2)(w_j\otimes v_i)\\
\phantom{\tPsi(z_1,z_2)}{}= q^{(\la|\nu)+(\la+2\rho|\wt(w_j)
+\wt(v_i))}(q^{-2\bar{\theta}(\la)}\otimes 1)q^{\pi_{W\otimes V}(\bar{T})}\\
\phantom{\tPsi(z_1,z_2)=}{}\times
R_{WV}(z_1/z_2)J_{WV}(z_1,z_2)(w_j\otimes v_i).
\end{gather*}
Combining 1) and 2), we obtain~\eqref{DiffK}.
\end{proof}

\subsection*{Acknowledgments}

The author would like to thank Michio Jimbo and Masato Okado for
stimulating discussions and valuable suggestions. He also thanks
Atsuo Kuniba and Atsushi Nakayashiki for discussions. He is also
grateful to the organizers of O'Raifeartaigh Symposium,  Janos
Balog, Laszlo Feher and  Zalan Horvath, for their kind invitation
and hospitality during his stay in Budapest.

\LastPageEnding

\end{document}